\documentclass{amsart}

\usepackage[english]{babel}

\usepackage[letterpaper,top=2cm,bottom=2cm,left=3cm,right=3cm,marginparwidth=1.75cm]{geometry}

\usepackage{enumitem}
\usepackage{amssymb}
\usepackage{amsmath}
\usepackage{amsthm}
\usepackage{graphicx}
\usepackage{physics}
\usepackage{tikz-cd}
\usepackage{mathtools}
\usepackage[all]{xy}
\usepackage{xcolor}
\usepackage{booktabs}
\usepackage{tabularx}
\usepackage{stmaryrd}

\usepackage{aliascnt}
\usepackage[colorlinks=true, allcolors=blue]{hyperref}
\usepackage[nameinlink]{cleveref}

\theoremstyle{plain}
\newtheorem{theorem}{Theorem}[section]
\newaliascnt{corollary}{theorem}
\newtheorem{corollary}[corollary]{Corollary}
\aliascntresetthe{corollary}
\newaliascnt{lemma}{theorem}
\newtheorem{lemma}[lemma]{Lemma}
\aliascntresetthe{lemma}
\newaliascnt{proposition}{theorem}
\newtheorem{proposition}[proposition]{Proposition}
\aliascntresetthe{proposition}
\newaliascnt{conjecture}{theorem}
\newtheorem{conjecture}[conjecture]{Conjecture}
\aliascntresetthe{conjecture}
\theoremstyle{definition}
\newaliascnt{definition}{theorem}

\aliascntresetthe{definition}
\newaliascnt{remark}{theorem}

\aliascntresetthe{remark}

\crefname{theorem}{Theorem}{Theorems}
\crefname{lemma}{Lemma}{Lemmas}
\crefname{proposition}{Proposition}{Propositions}
\crefname{corollary}{Corollary}{Corollaries}
\crefname{conjecture}{conjecture}{conjectures}
\Crefname{conjecture}{Conjecture}{Conjectures}
\crefname{assumption}{assumption}{assumptions}
\Crefname{assumption}{Assumption}{Assumptions}

\DeclareMathOperator{\Ad}{Ad}
\DeclareMathOperator{\Br}{Br}
\DeclareMathOperator{\End}{End}
\DeclareMathOperator{\Mat}{Mat}
\DeclareMathOperator{\Hom}{Hom}
\DeclareMathOperator{\Mon}{Mon}
\DeclareMathOperator{\SO}{SO}
\DeclareMathOperator{\OO}{O}

\DeclareMathOperator{\SU}{SU}

\newcommand{\NS}{\operatorname{NS}}

\newcommand{\Ext}{\operatorname{Ext}}

\newcommand{\cEnd}{\operatorname{\mathcal{E}nd}}

\newcommand{\C}{\mathbb C}
\newcommand{\PP}{\mathbb P}
\newcommand{\QQ}{\mathbb Q}
\newcommand{\R}{\mathbb R}
\newcommand{\Z}{\mathbb Z}

\newcommand{\cK}{\mathcal K}

\newcommand{\cM}{\mathcal M}

\newcommand{\cO}{\mathcal O}
\newcommand{\cP}{\mathcal P}

\newcommand{\cU}{\mathcal U}

\newcommand{\id}{\mathrm{id}}
\newcommand{\rk}{\operatorname{rk}}
\newcommand{\ch}{\operatorname{ch}}
\newcommand{\td}{\operatorname{td}}
\newcommand{\im}{\operatorname{im}}
\newcommand{\Sym}{\operatorname{Sym}}

\newcommand{\Pic}{\mathrm{Pic}}

\title[D-equivalence conjecture for hyperk\"ahler varieties of generalized Kummer type]{The D-equivalence conjecture for hyperk\"ahler varieties of generalized Kummer type}

\author[W. Donovan]{Will Donovan}
\address{Yau Mathematical Sciences Center,
Jingzhai, Tsinghua University, Haidian District,
Beijing, China, 100084}
\email{donovan@tsinghua.edu.cn}

\author[C. Wang]{Chengxi Wang}
\address{Yau Mathematical Sciences Center,
Jingzhai, Tsinghua University, Haidian District,
Beijing, China, 100084}
\email{chxwang@tsinghua.edu.cn}

\author[N. Zhang]{Nantao Zhang}
\address{N. Zhang: School of Mathematics, Sun Yat-sen University, Guangzhou 510275, China}
\email{zhangnt@mail.sysu.edu.cn}

\date{}

\begin{document}

\begin{abstract}
  We show that birationally equivalent hyperk\"ahler varieties of generalized Kummer type are derived equivalent, confirming the D-equivalence conjecture in this case.
  Combined with recent results of other authors, this settles the conjecture for all known deformation types of hyperk\"ahler varieties.
  The kernel is built from an equivariant Fourier-Mukai kernel on the sum-zero loci of powers of an abelian surface.
  We deform it along diagonal twistor paths, using the hyperholomorphic endomorphism algebra.
  Comparing wall hyperplanes gives a common twisted derived category for two adjacent birational models, and a finite sequence of wall crossings then gives the equivalence.
\end{abstract}

\maketitle

\setcounter{tocdepth}{1}
\tableofcontents

\section{Introduction}

For a smooth projective variety $X$, let \(D^b(X)\) be the derived category of bounded complexes of coherent sheaves on $X$. Exploring the behaviour of \(D^b(X)\) under birational transformations has attracted considerable interest over the years. Birational geometry and derived categories are two significant viewpoints to study smooth projective varieties. 
We recall the D-equivalence conjecture for algebraic varieties (\cite{BondalOrlov1995}, \cite[Conjecture 1.2]{Kawamata2002}) which proposes an underlying unity of the two viewpoints.
\begin{conjecture}
  Let \(X\) and \(X'\) be birational smooth projective varieties. Then \(D^b(X)\) and \(D^b(X')\) are equivalent if there is a common smooth projective resolution \(p: Z \to X\), \(p': Z \to X'\) with \(p^*K_X \sim (p')^*K_{X'}\).
\end{conjecture}
In particular, it implies that birational varieties with trivial canonical bundles are derived equivalent. The conjecture is known to hold for threefolds with trivial canonical bundle \cite{Bridgeland2002,Kawamata2002}, for hyperk\"ahler fourfolds \cite{Namikawa2003, WW2003}, for hyperk\"ahler varieties which can be realized as a Bridgeland moduli space of stable objects on a K3 surface \cite{Leistner2021}, for hyperk\"ahler varieties of \(K3^{ [n]}\)-type \cite{MSYZ2025}, and for hyperk\"ahler varieties of OG10 type \cite{HS2026}. In a recent preprint \cite{Ling2026}, it is claimed that birational hyperk\"ahler varieties of OG6 type are isomorphic. In this paper, we prove the conjecture for the remaining known families of hyperk\"ahler varieties, that is generalized Kummer type of arbitrary dimension.

\begin{theorem}
    Let \(X\) and \(X'\) be birational smooth projective complex hyperk\"ahler varieties of generalized Kummer type with complex dimension \(2n\), where \(n \geq 2\). Then there exists an equivalence
    \[D^b(X) \cong D^b(X').\]
\end{theorem}

Our method 
is inspired by Markman's work of deforming projectively hyperholomorphic bundles along generic diagonal twistor paths \cite{Markman2024, Markman2020}. Then we  take a similar strategy to \cite{MSYZ2025} of finding an intermediate birational model of two birational models with adjacent K\"ahler cones to bridge their twisted derived categories.

However, the construction of the Fourier-Mukai kernel in \cite{MSYZ2025} does not directly transfer to generalized Kummer type. We develop a new approach to obtain the initial equivariant equivalence from Yang's lifting and splitting results \cite[Theorem 2.5, Proposition 2.11]{Yang2025}. A key ingredient is Orlov's fundamental short exact sequence for derived equivalence of abelian varieties (\cite[Theorem~5.1]{Magni2026}, \cite[Section 4]{Orlov2002}). Then we choose an appropriate sign linearization so that the induced isometry matches K\"ahler classes, and deform the kernel projectively. Finally, the large parameter \(g\) makes the preimage of a chosen target wall fail to be a source wall. This gives the same twisted source for two adjacent birational models.

\subsection*{Acknowledgements}
The first author is grateful to Nicolas Addington for his insight and expertise over the years, and to Qizheng Yin for his explanations of \cite{MSYZ2025}.
The second author would thank Eyal Markman, Burt Totaro and Qizheng Yin for helpful conversations.
The third author is grateful to Shizhuo Zhang for his continuous support, and would like to thank Moritz Hartlieb, Saket Shah and Huishi Yu for insightful discussions and advice.

\subsection*{AI disclosure}
The project was started two years ago motivated by Markman's work of showing the algebraicity of rational Hodge isometries \cite{Markman2024,Markman2020} and the main methods in
\cite{MSYZ2025}. As the tools for finding the kernel in \cite{MSYZ2025} do not apply to generalized Kummer type, we came up with the idea of using lifting and splitting methods from Yang's work \cite{Yang2025}. With several rounds of help from ChatGPT 6.0,  we realized our idea to obtain a Fourier-Mukai kernel.
The authors independently verified the proof and take responsibility for the mathematical content. ChatGPT and Claude Opus~5.5 also helped fix typos and check the write-up of our manuscript.

\subsection*{Conventions}
We work over the complex numbers \(\C\). Let \(n \geq 2\) and \(m = n + 1 \geq 3\). If \(A\) is an abelian surface, we denote
\[N_A = \ker(A^m \xrightarrow{\sum} A) \quad\text{and}\quad K(A) = K_{n}(A) = (A^{[m]} \xrightarrow{\sum} A)^{-1}(0).\]
Then \(K(A)\) is a hyperk\"ahler variety of generalized Kummer type of dimension \(2n\). The symmetric group \(G = S_m\) acts on \(N_A\) by permutation. For a bundle of rank \(r\), we write
\[\kappa(E) = \ch(E) e^{-c_1(E)/r}\]
for the normalized Chern character and \(\kappa_i\) for the complex degree \(i\) part of \(\kappa\).

\section{Preliminaries: Abelian surfaces}

Let \(A\) be an abelian variety. That is a projective complex torus with a chosen origin \(0\). Its dual abelian variety is \(A^{\vee} = \Pic^{0}(A)\), with tensor product as the group law. We write \(P_{\xi}\) for the line bundle represented by \(\xi \in A^{\vee}\). The \emph{Poincar\'e bundle} \(\cP\) is the line bundle on \(A \times A^{\vee}\) such that \(\cP|_{A\times \{\xi\}} \cong P_{\xi}\). For \(x \in A\), write \(t_{x}(y) = y + x\). And for a line bundle $L$ on $A$, we define a homomorphism
\[\phi_{L}: A \to A^{\vee}, \quad x \to [t_{x}^{*}L \otimes L^{-1}].\]
Since \(\phi_{L}\) equals \(\phi_{L \otimes P}\) when \(P \in \Pic^{0}(A)\), we may write \(\phi_{H} \coloneq \phi_{L_{H}}\) with \(c_{1}(L_{H}) = H \in \operatorname{NS}(A)\). We denote \(A[n] := \{a \in A \mid na = 0\}\).

A vector bundle \(F\) on \(A\) is \emph{semihomogeneous} if for every \(x \in A\), there is \(P \in \Pic^{0}(A)\) with \(t_{x}^{*}F \cong F \otimes P\). It is called \emph{simple} if \(\End(F) = \C\). And its divisor slope is \(\mu(F) = c_{1}(F) / \rk(F) \in \NS(A)_{\QQ}\). For a fixed ample or K\"ahler class \(\omega\), the \emph{numerical slope} of \(F\) on an \(n\)-dimensional variety is \(\mu_{\omega}(F) = \int c_{1}(F) \omega^{n - 1}/\rk(F)\). A vector bundle \(F\) is called (semi)stable (with respect to \(\omega\)) if for every subsheaf \(E \subset F\) with \(0 < \rk(E) < \rk(F)\), we have
\[\mu_{\omega}(E) (\leq)< \mu_{\omega}(F).\]

A vector bundle \(F\) is called polystable if it is a direct sum of stable bundles of same slope.

We use simple semihomogeneous bundles to build locally free Fourier-Mukai kernels. The next result gives their existence and the rank.

\begin{theorem}[{\cite[Proposition 6.16, Theorem 7.11]{Mukai1978}, \cite[Proposition 4.6, Lemma 4.9]{Orlov2002}}]
  \label{thm:semihomo}
  Let \(\mu = c_{1}(L)/t \in \NS(A)_{\QQ}\), with \(t \in \Z_{> 0}\), and let
  \[\Phi_{\mu} = \im (A \to A \times A^{\vee}, x \to (tx, \phi_{L}(x))).\]
  There exists a simple semihomogeneous bundle \(F\) of slope \(\mu\). Every other such bundle is \(F \otimes P\) for some \(P \in \Pic^{0}(A)\). Write \(p_{1}: \Phi_{\mu} \to A\). One has
  \[\rk(F)^{2} = \deg p_{1}\]
  \[\{P: F \otimes P \cong F\} = \phi_{L}(A[t]).\]
  Moreover, such \(F\) is slope-stable for every ample polarization and \(\ch(F) = \rk(F) \exp(\mu)\).
\end{theorem}

We recall when such bundles define equivalences. Let \(A\) and \(B\) be two abelian varieties with an isomorphism
\[\gamma : A \times A^{\vee} \to B \times B^{\vee}.\]
We could write \(\gamma\) in matrix form \(
\begin{pmatrix}
  a & b \\
  c & d \\
\end{pmatrix}
\). Such isomorphism is called \emph{isometric} when
\[\gamma^{-1} =
  \begin{pmatrix}
    d^{\vee} &  -b^{\vee} \\
    -c^{\vee} & a^{\vee} \\
  \end{pmatrix}.
\]

\begin{theorem}[{\cite[Construction 4.10, Proposition 4.11-4.12]{Orlov2002}}]
  \label{thm:equivalence}
  Let \(A\), \(B\) be abelian varieties and \(\gamma: A \times A^{\vee} \to B \times B^{\vee}\) be an isometric isomorphism. Assume \(b : A^{\vee} \to B\) is an isogeny. The symmetric rational homomorphism
  \[
    \begin{pmatrix}
      b^{-1}a & -b^{-1} \\
      -(b^{-1})^{\vee} & db^{-1} \\
    \end{pmatrix}
  \]
  defines a class \(\mu \in \NS(A \times B)_{\QQ}\). Any simple semihomogeneous bundle \(P\) of this slope is a Fourier-Mukai equivalence kernel \(D^{b}(A) \to D^{b}(B)\), has \(\rk(P)^{2} = \deg b\) and induces \(\gamma\).
\end{theorem}

On a surface, the initial kernel is a universal family on a moduli space of stable sheaves. We first recall the numerical invariants of this moduli space.

For an abelian surface \(A\), we have Todd class \(\td_{A} = 1\). The \emph{Mukai vector} of a coherent sheaf \(F\) is
\[v(F) = \ch(F) \sqrt{\td_{A}} = (\rk(F), c_{1}(F), \chi(F)) \in H^{0}(A, \Z) \oplus H^{2}(A, \Z) \oplus H^{4}(A, \Z),\]
with the pairing
\[\langle (r, \xi, s), (r', \xi', s') \rangle = \int_{A} \xi \xi' - r s' - r' s.\]

A vector \(v\) with \(v^{2} = \langle v,v\rangle = 0\) is called isotropic. When the moduli space is non-empty, isotropy gives a surface, while coprimality of \(r\) and \(s\) gives a universal bundle.

\begin{proposition}
  \label{prop:surface-moduli}
  Let \(A\) be an abelian surface, \(L\) an ample line bundle, and \(v = (r, \xi, s) \in \Z \oplus \NS(A) \oplus \Z\) satisfy
  \[r > 0, \quad v^{2} = 0, \quad \gcd(r,s) = 1.\]
  Assume that the moduli space \(M = M_{A, L}(v)\) of slope stable sheaves of Mukai vector \(v\) is non-empty. Then \(M\) is a smooth projective abelian surface after choosing a point as origin, and every sheaf it parametrizes is a simple semihomogeneous vector bundle. There is a universal vector bundle \(\cU\) on \(M \times A\). Its Fourier-Mukai functor is an equivalence, and \(\cU\) is semihomogeneous on \(M \times A\).
\end{proposition}

\begin{proof}
  Stability implies simplicity. Since \(v^2 = 0\), \cite[Proposition 3.8]{Mukai1987Fourier} shows that every sheaf represented by \(M\) is a semihomogeneous vector bundle. Fix one such bundle \(F\). By \Cref{thm:semihomo}, every other one is \(F \otimes P\) for some \(P \in \Pic^0(A)\), and the subgroup \(\Sigma(F) = \{P \mid F \otimes P \cong F\}\) is finite. The trace map \(\Ext^2(F,F) \to H^2(A,\cO_A)\) is an isomorphism, so the deformations are unobstructed and \(M\) is smooth. Thus tensoring by \(\Pic^0(A)\) identifies \(M\) with \(\Pic^0(A)/\Sigma(F)\), which is a projective abelian surface after choosing an origin.

  The class of a point \(a \in A\) has a locally free resolution. So the gcd of \(\chi(F \otimes V)\), with \(V\) ranging over vector bundles on \(A\), divides both \(\chi(F) = s\) and \(\chi(F \otimes \cO_a) = r\). Since \(\gcd(r,s) = 1\), \cite[Theorem A.6]{Mukai1987Moduli} gives a universal vector bundle \(\cU\). Its distinct fibers are orthogonal by \cite[Lemma 4.8]{Orlov2002}, and their self-Hom spaces are one-dimensional. So \cite[Theorem 1.5]{Orlov2002} gives full faithfulness, and the trivial canonical bundle gives the equivalence. Now \cite[Theorem 2.10, Corollary 2.13]{Orlov2002} identifies the translation and tensor stabilizer of \(\cU\), up to a sign on the source dual factor, with the graph of the associated isomorphism. This stabilizer has dimension \(4\), and its projection to \(M \times A\) has finite kernel since \(\cU\) is locally free. So it surjects onto \(M \times A\), which proves semihomogeneity on the product.
\end{proof}

We record the McKay correspondence on the zero fiber, with the finite flatness and duality needed for the kernel, the determinant slopes, and the sign calculation.

\begin{proposition}
  \label{prop:kummer-mckay}
  Let \(A\) be an abelian surface, \(G=S_m\), and let \(\Gamma_A\) be the restriction to the summation fiber over \(0\) of the reduced isospectral Hilbert-Chow correspondence. The projections
  \[N_A \xleftarrow{b_A} \Gamma_A \xrightarrow{h_A} K(A)\]
  satisfy that \(h_A\) is finite flat of degree \(m!\) and \(\Gamma_A\) is Gorenstein. Its \(G\)-linearized relative dualizing line is
  \[\omega_{h_A}\cong {\mathrm{sgn}}\otimes h_A^*\cO_{K(A)}(\delta_A).
  \]
  The functor
  \[{\mathrm{BKR}}_A=(h_{A*}Lb_A^*)^G:D_G^b(N_A)\xrightarrow{\sim}D^b(K(A))
  \]
  is an equivalence, with inverse
  \[{\mathrm{BKR}}_A^{-1}(F)=Rb_{A*}(\omega_{h_A}\otimes h_A^*F).
  \]
\end{proposition}

\begin{proof}
  Let \(Z_A=(A^m\times_{A^{(m)}}A^{[m]})_{\mathrm{red}}\) be the isospectral Hilbert scheme, with projections \(b\) and \(h\). Write \(j:N_A\hookrightarrow A^m\) and \(i:K(A)\hookrightarrow A^{[m]}\) for the zero-fiber inclusions. The restricted correspondence fits into the commutative diagram
  \[
    \begin{tikzcd}[column sep=large, row sep=large]
      N_A \arrow[d, hook, "j"'] & \Gamma_A \arrow[l, "b_A"'] \arrow[r, "h_A"] \arrow[d, hook] & K(A) \arrow[d, hook, "i"] \\
      A^m & Z_A \arrow[l, "b"'] \arrow[r, "h"] & A^{[m]}.
    \end{tikzcd}
  \]
  Both squares are Cartesian, since the top row is the fiber over \(0\) of the bottom row under summation. Translation identifies the pullback of the bottom correspondence along the finite \'etale map \([m]:A\to A\) with
  \[A\times N_A\xleftarrow{\id_A\times b_A}A\times\Gamma_A
    \xrightarrow{\id_A\times h_A}A\times K(A).
  \]
  Hence the finite-flat and Gorenstein properties and the relative-duality formula descend from \cite[Theorem 3.1, Proposition 3.7.4]{Haiman2001}. The sign records the alternating character of Haiman's perfect pairing.

  For the Hilbert-Chow McKay correspondence \cite[Theorem 1.1, Corollary 1.3]{BKR2001}, we use Krug's transform
  \[\Psi^{[m]}=(h_*Lb^*)^G:D_G^b(A^m)\longrightarrow D^b(A^{[m]}).
  \]
  This is an equivalence whose restriction to the zero fiber is \({\mathrm{BKR}}_A\), by \cite[Proposition 2.8 and Remark 3.14]{Krug2018}, see also \cite[Lemma 6.2]{Meachan2015}. In particular, the diagram gives the compatibility
  \[Li^*\circ\Psi^{[m]}\cong {\mathrm{BKR}}_A\circ Lj^*.
  \]
  Finite duality identifies the right adjoint of \({\mathrm{BKR}}_A\), and hence its inverse, with \(Rb_{A*}h_A^!\), where \(h_A^!=\omega_{h_A}\otimes h_A^*\).
\end{proof}

\section{Preliminaries: Hyperk\"ahler geometry}
\label{sec:prel-hyperk-geom}

We recall some basic theory of hyperk\"ahler geometry, see \cite{Debarre2022}, for example. An \emph{irreducible holomorphic symplectic manifold} (IHSM in short) is a simply connected compact K\"ahler manifold \(X\) with \(H^0(X,\Omega_X^2) = \C\sigma_X\), where \(\sigma_X\) is an everywhere nondegenerate holomorphic two-form. We call a smooth projective variety whose underlying complex manifold is an IHSM a \emph{hyperk\"ahler variety}. IHSMs need not be projective, and we also consider them along twistor paths.

Let \(A\) be an abelian surface and \(m = n + 1\). The Hilbert scheme \(A^{[m]}\) parametrizes zero-dimensional subschemes of length \(m\). The Hilbert-Chow morphism sends such a subscheme to its associated zero-cycle in \(A^{(m)} = A^m/S_m\), counting points with their multiplicities. Composing it with addition gives
\[s_m: A^{[m]} \xrightarrow{\mathrm{HC}} A^{(m)} \xrightarrow{\sum} A.\]
The \emph{generalized Kummer variety} is \(K_n(A) = s_m^{-1}(0)\). It is a hyperk\"ahler variety of dimension \(2m - 2 = 2n\). The restriction of the Hilbert-Chow morphism is \(K(A) \to N_A/S_m\), which is the quotient used in \Cref{prop:kummer-mckay}. An IHSM is \emph{of generalized Kummer type} if it is deformation equivalent to a \(K_n(A)\).

For an IHSM \(X\) of dimension \(2n\), we have a canonical nondegenerate symmetric bilinear primitive pairing called \emph{Beauville-Bogomolov-Fujiki form} (BBF form in short) on \(H^{2}(X, \Z)\),
\[(,)_{X}:H^{2}(X, \Z) \times H^{2}(X, \Z) \to \Z.\]
We write \(q_{X}(u) = (u, u)_{X}\) for the pairing. The pairing is characterized by the Fujiki relation
\[\int_{X} u^{2n} = c_{X}q_{X}(u)^{n},\]
where \(c_{X} > 0\) is the constant.

The BBF form also describes the Lie algebra acting on cohomology. Let \(\mathfrak g_{X}\) be its \emph{Looijenga-Lunts-Verbitsky Lie algebra} (LLV algebra in short) \cite{LL1997,Verbitsky1996coh}. We define
\[\tilde{H}(X, \QQ) = \QQ \alpha \oplus H^{2}(X, \QQ) \oplus \QQ \beta, \quad (\alpha, \beta) = -1, \alpha^{2} = \beta^{2} = 0.\]
Then the LLV structure theorem (\cite[Proposition 4.5]{LL1997}; \cite[Theorem 1.2]{Verbitsky1996coh}) identifies two Lie algebras
\[\mathfrak g_{X} \cong \mathfrak{so}(\tilde{H}(X, \QQ)).\]
For an operator of degree \(i\), one sends \(H^{j}\) to \(H^{j + i}\). The structure theorem gives
\[\mathfrak{g}_{X} = (\mathfrak g_{X})_{-2} \oplus (\mathfrak g_{X})_{0} \oplus (\mathfrak g_{X})_{2}, \quad (\mathfrak g_{X})_{0} = \QQ h_{X} \oplus \bar{\mathfrak{g}}_{X}\]
where
\[\bar{\mathfrak{g}}_{X} = \mathfrak{so}(H^{2}(X, \QQ)).\]

For the generalized Kummer variety, \(H^{2}(K(A),\Z)\) is given by two orthogonal summands
\[H^{2}(K(A), \Z) = H^{2}(A, \Z) \oplus \Z \delta_{A}\]
where \(H^{2}(A, \Z) \cong U^{\oplus 3}\) and \(\delta_{A}^{2} = -2m\). \(2 \delta\) is the Hilbert-Chow exceptional divisor.

Now let \(SH^{*}(X, \QQ)\) be the subalgebra generated by \(H^{2}(X, \QQ)\). It is an irreducible LLV subrepresentation \cite[\S1]{Verbitsky1996coh}. For a fixed complex structure, define its \emph{infinitesimal Hodge operator} \(W_{X}\) by
\[W_{X}|_{H^{p,q}(X)} = i(p - q) {\mathrm{id}}.\]
It belongs to \(\bar{\mathfrak{g}}_{X} \otimes \R\). In particular, in degree \(2k\), its kernel is exactly \(H^{k, k}(X)\). For an invertible map \(u\), we write \(\Ad_{u}(\xi) = u\xi u^{-1}\).

These structures let us study Fourier-Mukai equivalences through their action on cohomology. We recall the compatibility with the LLV algebra:

\begin{theorem}[{\cite[Theorems A--C, 4.10--4.11]{Taelman2023}}]
  \label{thm:LLV}
  Let \(X,Y\) be projective varieties of generalized Kummer type of dimension \(2n\), and let \(\Phi:D^b(X)\xrightarrow{\sim}D^b(Y)\) be a Fourier-Mukai equivalence. Its cohomological transform \(\Phi^H\) conjugates their rational LLV algebras and restricts to a Mukai isometry \(\Phi^{SH}:SH^*(X,\QQ)\xrightarrow{\sim}SH^*(Y,\QQ)\).

  Under the LLV identifications \(SH^*\cong\Sym_0^n\tilde H\), normalized by \(1\mapsto\alpha^n/n!\), there are a rational Hodge isometry \(g_\Phi:\tilde H(X,\QQ)\xrightarrow{\sim}\tilde H(Y,\QQ)\) and a sign \(a_\Phi\in\{\pm1\}\) such that
  \[\Phi^{SH}=a_\Phi\Sym_0^n(g_\Phi).\]
  Here \(\Sym_0^n\) is the trace-free symmetric power and \(\alpha,\beta\) have Hodge type \((1,1)\). The projective class \([g_\Phi]\), modulo \(g_\Phi\sim-g_\Phi\), is uniquely determined and functorial in \(\Phi\). Replacing \(g_\Phi\) by \(-g_\Phi\) replaces \(a_\Phi\) by \((-1)^n a_\Phi\).
\end{theorem}

For kernels of nonzero rank, the normalized Chern character isolates the isometry on \(H^{2}\). We use the following form of this result.

\begin{lemma}[{\cite[Proposition 3.1, Lemmas 3.7, 4.1]{Markman2024}}]
  \label{lem:degree-reversing}
  Let \(E\) be a nonzero-rank Fourier-Mukai equivalence kernel between projective varieties \(X,Y\) of generalized Kummer type of dimension \(2n\). The normalized cohomological transform
  \[\phi_{E} = [\kappa(E)\sqrt{\td_{X \times Y}}]_{*}: H^{*}(X, \QQ) \to H^{*}(Y, \QQ)\]
  is degree reversing. Choose a representative \(g_{\phi_E}\) of its projective extended isometry as in \Cref{thm:LLV}. Its restriction to \(H^2\) is a rational Hodge isometry. We write \(\psi_E\) for the unique choice among \(\pm g_{\phi_E}|_{H^2}\) that maps the component of \(\{x\in H^{1,1}(X,\R)\mid q_X(x)>0\}\) containing the K\"ahler cone to the corresponding component for \(Y\). Then
  \[{\mathrm{Ad}}_{\phi_{E}}|_{\bar{\mathfrak g}_{X}} = {\mathrm{Ad}}_{\psi_{E}}.\]

  Moreover, \(\psi_{E}\) is a nonzero scalar multiple of
  \[H^{2}(X, \QQ) \xrightarrow{\cup c_{2}(X)^{n - 1}} H^{4n - 2}(X, \QQ) \xrightarrow{\phi_{E}} H^{2}(Y, \QQ)\]
\end{lemma}

To deform this isometry together with a kernel, we work with marked pairs of manifolds. A \emph{marking} is an integral isometry \(\eta: H^{2}(X, \Z) \to \Lambda\). Its period is the complex line \([\eta(\sigma_{X})]\), where \(\sigma_{X}\) is a nonzero holomorphic symplectic form. The period domain is
\[\Omega_{\Lambda} = \{ [z] \in \PP(\Lambda_{\C}) \mid (z, z) = 0, (z, \bar{z}) > 0\}.\]
A \emph{parallel-transport operator} is the cohomology identification obtained from the Gauss-Manin local system along a path in a smooth proper family, or a composition of such identifications. Fix a rational lattice isometry \(\psi \in \OO(\Lambda_{\QQ})\). The space \(\cM_{\psi}\) consists of marked quadruples \(q = (X_{1}, \eta_{1}, X_{2}, \eta_{2})\) such that
\[\psi_{q} = \eta_{2}^{-1}\psi \eta_{1}: H^{2}(X_{1}, \QQ) \to H^{2}(X_{2}, \QQ)\]
is a Hodge isometry and some \(\omega \in \cK_{X_{1}}\) satisfies \(\psi_{q}(\omega) \in \cK_{X_{2}}\). The last condition is called \emph{K\"ahler compatibility}. A superscript \(0\) denotes one specified connected component.

A hyperk\"ahler metric has complex structures \(I,J,K\) satisfying the quaternion relations. Its \emph{twistor sphere} consists of the complex structures \(aI + bJ + cK\) with \(a^{2} + b^{2} + c^{2} = 1\) on the same differentiable manifold. A \emph{diagonal twistor line} rotates both factors of a compatible pair simultaneously, with metrics chosen from the matched K\"ahler classes. A \emph{twistor path} is a finite chain of such lines. A \emph{joint} is a point shared by consecutive lines. The path is \emph{generic} when its interior joints have trivial Picard groups.

The next result lets us choose compatible pairs and join them by such paths.

\begin{lemma}[{\cite[Lemma 5.7--5.8]{Markman2024}}]
  \label{lem:pair-paths}
  For a connected component \(\cM_{\psi}^{0}\), each projection onto the corresponding component of marked IHSM moduli is surjective. Every two points of \(\cM_{\psi}^{0}\) can be joined by a generic diagonal twistor path.
\end{lemma}

We now recall how to continue bundles along these paths. The unit quaternion group \(\SU(2)\) acts on differential forms of a hyperk\"ahler manifold. For a two-form, invariance is equivalent to having type \((1, 1)\) for every complex structure on the twistor sphere. A connection is \emph{hyperholomorphic} if its curvature has this property. It is \emph{projectively hyperholomorphic} when its projective curvature, or equivalently the traceless curvature of a local lift, has this property. The above two notions were introduced in \cite{Verbitsky1996}.

An \emph{Azumaya algebra} of degree \(R\) is a sheaf of algebras locally isomorphic to \(\Mat_{R}(\cO)\). It determines an associated projective-space bundle. A local splitting is an identification with \(\cEnd(F)\) for a rank \(R\)-bundle \(F\). These local splittings define a possibly twisted bundle whose endomorphism algebra is the given Azumaya algebra.

We use the following continuation result for the endomorphism algebra of a kernel on a compatible pair. The algebra continues, while its holomorphic splittings may become twisted.

\begin{proposition}[{\cite[Proposition 5.15]{Markman2024}}]
  \label{prop:azumaya-continue}
  Let \(q = (X_{1}, \eta_{1}, X_{2}, \eta_{2}) \in \cM^{0}_{\psi}\) and let \(F\) be a locally free sheaf on \(X_{1} \times X_{2}\). Assume that:
  \begin{enumerate}
    \item \(F\) is stable for \({\mathrm{pr}}_{1}^{*} \omega + {\mathrm{pr}}_{2}^{*}\psi_{q}(\omega)\) on a non-empty open subcone of matched K\"ahler classes.
    \item the transported \(c_{2}(\cEnd F)\) is of type \((2, 2)\) along the compatible deformations under consideration.
  \end{enumerate}
  Along a generic diagonal twistor path, whose first line uses this stable cone, \(\cEnd F\) continues as a locally free Azumaya algebra. Each endpoint has the form \(\cEnd F'\) for a possibly twisted locally free bundle \(F'\). If condition (2) holds throughout the compatible component, every point of the component can be reached in this way.
\end{proposition}

For the kernel constructed below, we invoke continuation as a twisted Fourier-Mukai equivalence directly in \Cref{prop:ext-continue}. We first fix the notation for the Brauer class of a rational degree-two class. For \(B \in H^{2}(X, \QQ)\), we denote by \([B]_{X}\) its image under
\[H^{2}(X, \QQ) \to H^{2}(X, \cO_{X}) \to H^{2}(X, \cO_{X}^{*})\]
where the second arrow is induced by the exponential exact sequence. On a compact K\"ahler manifold, the first arrow is projection to \(H^{0, 2}\). The endpoint Brauer classes will be computed from the transported determinant slopes in \Cref{sec:wall-crossing}, using \cite[Propositions 4.2--4.3]{Caldararu2002}.

The final argument also needs control of markings and K\"ahler chambers. We recall the monodromy group and the relation between periods and birational models.

For an even lattice \(L\), its dual lattice is
\[L^{*} = \{v \in L \otimes \QQ \mid (v, L) \subset \Z\},\]
with the discriminant group
\[A_{L} = L^{*}/L.\]
The finite group \(A_{L}\) has the induced quadratic form with values in \(\QQ/2\Z\), and \(\OO(A_{L})\) denotes its automorphism group preserving that form. Let \(\widetilde{\SO}^{+}(L) = \ker(\SO^{+}(L) \to \OO(A_{L}))\). The group \(\Mon^{2}(X) \subset \OO(H^{2}(X, \Z))\) is generated by degree \(2\) parallel transport around loops in smooth proper families based at \(X\).

For generalized Kummer type, the monodromy group has a lattice theoretic description.

\begin{theorem}[{\cite[Theorem 4.3]{Mongardi2016}}]
  \label{thm:kummer-monodromy}
  Let \(X\) be an IHSM of generalized Kummer type of dimension \(2n\), \(n \geq 2\), with \(H^{2}(X, \Z) \cong \Lambda\). Let
  \[W(\Lambda) = \{u \in \OO^{+}(\Lambda) \mid u|_{A_{\Lambda}} = \pm 1\},\]
  and \(\epsilon_{\mathrm{disc}}(u)\) the sign of \(u|_{A_{\Lambda}}\). Then \(\Mon^{2}(X) = \ker (\det \cdot \epsilon_{\mathrm{disc}}: W(\Lambda) \to \{\pm 1\})\). In particular, \(\widetilde{\SO}^{+}(\Lambda) \subset \Mon^{2}(X)\).
\end{theorem}

Recall that two points are \emph{inseparable} if they have no disjoint open neighborhoods.
The next result describes the period map and the markings of birational models.

\begin{theorem}[{\cite[Theorems 2.2, 3.1]{Markman2011}}]
  \label{thm:torelli}
  On each connected component of marked IHSM moduli the period map is surjective, and its Hausdorff reduction is identified with the period domain. Points of that component with the same period are inseparable. The degree-two cohomology identification induced by a birational map between IHSMs is a parallel-transport operator. With markings corresponding under this identification, the two marked manifolds lie in the same connected component and are inseparable.
\end{theorem}

The following projectivity criterion will be used for the common source.

\begin{theorem}[{\cite[Theorem 2]{Huybrechts2003}}]
  \label{thm:projectivity}
  An IHSM is projective if and only if it has an integral \((1, 1)\) class of positive BBF square.
\end{theorem}

We organize the chamber results as in \cite{MSYZ2025}, to make the reduction to adjacent birational models explicit. Let \(\mathcal C_X^+\) be the component of \(\{x\in H^{1,1}(X,\R):q_X(x)>0\}\) containing the K\"ahler cone \(\cK_X\). A \emph{wall} is a codimension-one boundary face of a K\"ahler chamber inside \(\mathcal C_X^+\). Its supporting hyperplane is \(e^{\perp}\) for a negative class \(e\). A negative integral class is \emph{monodromy birational minimal} (MBM in short) if, up to monodromy, its perpendicular supports a K\"ahler wall of a birational model. A \emph{K\"ahler type chamber} is a K\"ahler cone of a birational model transported to \(H^{1,1}(X, \R)\) by a birational identification and a Hodge monodromy operator.

The next result identifies the chamber walls and their behaviour under deformation.

\begin{theorem}[{\cite[Theorems 1.17, 1.19]{AV2015}, \cite[Theorem 5.16]{Markman2011}}]
  \label{thm:mbm-chambers}
  For an IHSM \(X\), let \(\mathcal S_{X}\) be the set of primitive integral MBM classes of type \((1, 1)\). The K\"ahler cone is a connected component of
  \[\mathcal C_{X}^{+}\setminus \bigcup_{e \in \mathcal S_{X}}e^{\perp}.\]
  Every component is a K\"ahler type chamber: it is the image of \(\cK_Y\) for an IHSM \(Y\) birational to \(X\) under a Hodge parallel-transport isometry \(H^2(Y,\Z)\to H^2(X,\Z)\). With the corresponding markings, \(X\) and \(Y\) have the same period and lie in the same connected component of marked moduli. The MBM property is invariant under deformations on which the class stays of type \((1, 1)\).
\end{theorem}

For a projective IHSM \(X\), define its \emph{birational K\"ahler cone} and the interior of its closure by
\[\mathcal{BK}_X=\bigcup_{f:X\dashrightarrow Y}f^*\cK_Y,\qquad
  \mathcal B_X={\mathrm{int}}_{\mathcal C_X^+}\!\left(\overline{\mathcal{BK}_X}^{\,\mathcal C_X^+}\right),\]
where \(Y\) ranges over IHSMs birational to \(X\), and the closure is taken inside \(\mathcal C_X^+\). The next result identifies the chambers coming from birational maps, together with their induced cohomology identifications.

\begin{theorem}[{\cite[Proposition 5.6, Corollary 5.7, Theorem 6.17]{Markman2011}}]
  \label{thm:birational-chambers}
  The cone \(\mathcal B_X\) is open and convex, and contains \(\mathcal{BK}_X\) as a dense subset. The connected components of
  \[\mathcal B_X\setminus\bigcup_{e\in\mathcal S_X}e^{\perp}\]
  are precisely the cones \(f^*\cK_Y\) for birational maps \(f:X\dashrightarrow Y\). These birational models \(Y\) are projective.
\end{theorem}

In Markman's terminology, \(\mathcal B_X\) is the fundamental exceptional chamber.

The squares of primitive MBM classes are also uniformly bounded. We use this bound in the lattice argument to exclude source wall classes.

\begin{theorem}[{\cite[Corollary 1.4]{AV2020}}]
  \label{thm:bound-q}
  For a fixed IHSM deformation type with \(b_{2} \geq 5\), there is a constant \(B > 0\) such that every primitive integral MBM class \(e\) satisfies
  \[0 < -q(e) \leq B.\]
  The constant depends on the deformation type, not on its period or Picard lattice.
\end{theorem}

In particular, such a \(B_{m}\) exists for IHSMs of generalized Kummer type of dimension \(2n\).

Local finiteness follows from \Cref{thm:bound-q} and the lattice argument of \cite[Remark 8.2.3]{Huybrechts2016}, as also recalled in \cite{MSYZ2025}.

\begin{lemma}
  \label{lem:mbm-local-finiteness}
  For an IHSM \(X\) of generalized Kummer type, the distinct hyperplanes \(e^{\perp}\), \(e\in\mathcal S_X\), are locally finite inside \(\mathcal C_X^+\). Consequently, the chamber decomposition is locally polyhedral there.
\end{lemma}

Together with \Cref{thm:birational-chambers}, this gives the standard reduction to adjacent birational models used in \cite{MSYZ2025}. See also the argument with paths between chambers in \cite[Remark 8.2.5(i)]{Huybrechts2016}.

\begin{corollary}
  \label{cor:adjacent-birational-models}
  Let \(X,X'\) be birational projective IHSMs of generalized Kummer type. Either they are isomorphic, or there is a finite sequence of birational projective models
  \[X=X_0,X_1,\ldots,X_s=X'\]
  such that, after identifying cohomology by birational maps, consecutive K\"ahler cones share a wall of codimension \(1\), supported by \(w_i^{\perp}\) for a primitive algebraic MBM class \(w_i\). Each wall can be crossed at a point lying on no other distinct MBM hyperplane.
\end{corollary}

Finally, we recall how to replace a marked factor by an inseparable one while staying in the same component of compatible pairs.

\begin{lemma}[{\cite[Lemma 2.2]{MSYZ2025}}]
  \label{lem:inseparable-pairs}
  Suppose \((X, \eta_{X}, Y, \eta_{Y}) \in \cM_{\psi}^{0}\), and \((X, \eta_{X})\) and \((X', \eta_{X'})\) are inseparable marked points. If \((X', \eta_{X'}, Y, \eta_{Y})\) is K\"ahler compatible with \(\psi\), it lies in the same component \(\cM_{\psi}^{0}\).
\end{lemma}

\section{Construction of equivariant kernel}
\label{sec:construct}

We use Yang's results on the equivariant lift and splitting \cite[Theorem 2.5, Proposition 2.11]{Yang2025}, and record the vector bundle normalization and numerical data needed for the deformation.

Fix a positive integer \(g\), set
\[c = 4m^2,\quad r = c^2g,\quad s = 1 + \frac{mr}{4},\quad d = gs.\]
Choose a Picard rank \(1\) abelian surface \(A\) with primitive type-\((1, d)\) polarization \(H\). Then \(H^2 = 2gs\) and \(v = (r, cH, s)\) satisfies \(v^2 = 0\). Notice that \(r\) and \(s\) are coprime. \(v\) is primitive. Set \(t = r/c = cg = 4m^2g\), so \(cH/r = H/t\) and \(\gcd(s,t) = 1\). By \Cref{thm:semihomo}, there is a simple semihomogeneous bundle of this slope. In Smith coordinates, \(\phi_H(A[t])\) has \(t^2(t/g)^2 = r^2\) elements, so the bundle has rank \(r\). Its Chern character is \(r\exp(H/t)\), and its Euler characteristic is \(rgs/t^2 = s\). Thus the moduli is non-empty.

Let \(M\) be the corresponding moduli surface and \(\cU\) the universal bundle on \(M \times A\), as in \Cref{prop:surface-moduli}. For the projections \(p_M,p_A\) from \(M\times A\), its Fourier-Mukai equivalence is
\[\Phi_{\cU}:D^b(M)\xrightarrow{\sim}D^b(A),\qquad
  \Phi_{\cU}(F)=\mathbf R(p_A)_*(p_M^*F\otimes\cU).\]
All stable sheaves with this Mukai vector are locally free. Otherwise, its reflexive hull would have negative Mukai square \(v(F^{**})^2 = v(F)^2 - 2r {\mathrm{length}}(F^{**}/F) < 0\), which contradicts Bogomolov's inequality.

The induced integral algebraic Mukai isometry identifies \(U \oplus \NS(M)\) with \(U \oplus \NS(A)\). Thus \(\NS(M) = \Z \hat{H}\) for a primitive ample class \(\hat{H}\) with \(\hat{H}^2 = H^2\). On the product, \(\cU\) is simple semihomogeneous, and we write
\begin{equation}
    \label{eq:chi}
    \frac{c_1(\cU)}{r} = a_M + a_A + \chi
\end{equation}
where
\[a_M = \frac{k}{r}\hat{H}, \quad a_A = \frac{c}{r}H, \quad \chi \in H^1(M, \QQ) \otimes H^1(A, \QQ),\quad k,c \in \Z.\]

Let \(f=\Phi_{\cU}^{H}:H^*(M,\QQ)\xrightarrow{\sim}H^*(A,\QQ)\) be its cohomological transform. Since both surfaces have Todd class \(1\), it is given by
\[f(x)=(p_A)_*(p_M^*x\cup\ch(\cU)).\]
Write \(\mathsf{e}_a(x)=e^a\cup x\) for the exponential cup product operator. Define the normalized transform
\[f_0\coloneq\mathsf{e}_{-a_A}\circ f\circ\mathsf{e}_{-a_M}=[re^{\chi}]_*.\]
The equality follows from \(\ch(\cU)=r\exp(a_M+a_A+\chi)\). The map \(f_0\) sends \(H^j(M,\QQ)\) to \(H^{4-j}(A,\QQ)\). Since \(f\) and the exponential cup product operators preserve the Mukai pairing, its restriction to \(H^2\) is a rational Hodge isometry. Define
\begin{equation}
  \label{eq:surface-isometry}
  \psi_{\cU}\coloneq\varepsilon_{\cU}f_0|_{H^2(M,\QQ)}:
  H^2(M,\QQ)\xrightarrow{\sim}H^2(A,\QQ),
\end{equation}
where \(\varepsilon_{\cU}\in\{\pm1\}\) is the unique sign for which \(\psi_{\cU}(\hat H)=H\). Equivalently, for \(x\in H^2(M,\QQ)\),
\[\psi_{\cU}(x)=\frac{\varepsilon_{\cU}r}{2}(p_A)_*(p_M^*x\cup\chi^2).\]

The equivalence \(\Phi_{\cU}\) also determines an isomorphism
\[\gamma_{\cU} = \begin{pmatrix}
    x & y \\
    z & w \\
  \end{pmatrix}: M \times M^{\vee} \xrightarrow{\sim} A \times A^{\vee}\]
\cite[Theorem 2.10 and Corollary 2.13]{Orlov2002}. This is an isometric isomorphism by \cite[Proposition 2.18]{Orlov2002}, and
\[\gamma_{\cU}^{-1} =
  \begin{pmatrix}
    w^{\vee} & -y^{\vee} \\
    -z^{\vee} & x^{\vee} \\
  \end{pmatrix}
\]
For an abelian variety \(A\) and a point \((a, \xi) \in A \times A^{\vee}\), we write
\[T_{(a, \xi)}^{A}(F) = t_{a *}F \otimes P_{\xi} = t_{-a}^{*}F \otimes P_{\xi}.\]
Then \(\gamma_{\cU}\) satisfies that
\[\Phi_{\cU} \circ T^{M}_{(a, \xi)} \cong T^{A}_{\gamma_{\cU}(a, \xi)} \circ \Phi_{\cU}.\]
For \(a\in M\), set \(\cU_a=\cU|_{\{a\}\times A}\), so that \(\Phi_{\cU}(\cO_a)\cong\cU_a\). The compatibility above gives
\[\Phi_{\cU}^{-1}(\cU_{a} \otimes P_{\beta}) \cong T^{M}_{(-y^{\vee}\beta, x^{\vee}\beta)}(\cO_{a}) \cong \cO_{a - y^{\vee} \beta}.\]
Consequently,
\begin{equation}
  \label{eq:ker-b}
  \ker y^{\vee} = \{\beta \in A^{\vee}: \cU_{a} \otimes P_{\beta} \cong \cU_{a}\} = \phi_{H}(A[t]).
\end{equation}
The last equality follows from \Cref{thm:semihomo}, and it has \(r^{2}\) elements. Since \(A^{\vee}\) and \(M\) have the same dimension, we obtain that the morphism \(y: M^{\vee} \to A\) is an isogeny with degree \(r^2\).

\begin{lemma}
  \label{lem:b-m2}
  The homomorphism \(y\) is divisible by \(m^{2}\) in \(\Hom(M^{\vee}, A)\).
\end{lemma}

\begin{proof}
  In the fixed target Smith coordinates, equation \eqref{eq:ker-b} becomes
  \[\ker y^{\vee} = (\Z/t\Z)^{2} \oplus (g (\Z/t\Z))^{2}\subset A^{\vee}[t] \cong (\Z/t\Z)^{4}.\]
  Since \(t = 4m^{2}g\), the inclusion \(A^{\vee}[m^{2}] \subset A^{\vee}[t]\) is contained in \(\ker y^{\vee}\). Therefore \(y^{\vee}\) factors through \([m^{2}]_{A^{\vee}}\). Dualizing proves the lemma.
\end{proof}

Thus \(\gamma_{\cU}\) satisfies the off-diagonal divisibility condition in Yang's lifting criterion \cite[Lemma 2.8]{Yang2025}, with his \(n\) equal to our \(m\).

We next consider the finite \'etale isogeny
\[q_A: N_A \times A \to A^m, \quad ((x_i), z) \to (x_i + z)_{i = 1}^{m}.\]
The kernel is \(((-z, \dots, -z), z)\) with \(mz = 0\). So its degree is \(m^4\). The \(S_m\)-action that is trivial on the second factor \(A\) fixes every point of the kernel. Now identify \(A^{n}\) with \(N_{A}\) by
\[(y_{1}, \dots, y_{n}) \to (y_{1}, \dots, y_{n}, - \sum_{i = 1}^{n} y_{i}),\]
and let \(e_{1}, \dots, e_{m}\) be the standard basis of \(\Z^{m}\). The integral matrix \(T\) of \(q_{A}\) has columns \(e_{1} - e_{m}, \dots, e_{n} - e_{m}\) followed by the diagonal column \((1, \dots, 1)^{t}\). Its determinant is \(m\). The first \(n\) columns have pairings \(\delta_{ij} + 1\). Let \(I_{n}\) denote the identity matrix and \(J_{n}\) the matrix with every entry equal to \(1\). Then
\[T^{t}T = {\mathrm{diag}}(G_{n}, m), G_{n} = I_{n} + J_{n}.\]

Set \(F=\cU^{\boxtimes m}\) on \(M^m\times A^m\). The next proposition applies Yang's lifting theorem and records the slope and rank of the resulting vector bundle.

\begin{proposition}
  \label{prop:rk-p}
    There is a diagonally \(S_m\)-linearized simple semihomogeneous equivalence kernel \(P\) on \(N_M \times N_A\), which is locally free, with
    \[\rk P = p = \frac{r^n}{m^2}.\]
    Its slope is the sum-zero part of the pullback slope of \(F\).
\end{proposition}

\begin{proof}
  By \Cref{lem:b-m2}, \(y=m^2y_0\) for a homomorphism \(y_0:M^{\vee}\to A\). Thus \(\Phi_{\cU}\) satisfies Yang's lifting criterion \cite[Lemma 2.8]{Yang2025}. The construction and splitting in \cite[Theorem 2.5 and its proof]{Yang2025} give a Fourier-Mukai equivalence \(D^b(N_M)\to D^b(N_A)\) with a compatible \(S_m\)-equivariant structure. Let \(\mathcal Q\in D^b(N_M\times N_A)\) be its diagonally \(S_m\)-linearized kernel.

  We first write its isometric map in our coordinates. The external product kernel \(F=\cU^{\boxtimes m}\) has associated map
  \[\gamma_F=\begin{pmatrix}xI_m&yI_m\\zI_m&wI_m\end{pmatrix}.\]
  The isogenies \(q_M,q_A\) are represented by \(T\), and their duals by \(T^t\). Transporting \(\gamma_F\) to the product coordinates therefore gives upper-right and lower-left blocks
  \[T^{-1}(yI_m)(T^t)^{-1}=y\,{\mathrm{diag}}(G_n^{-1},m^{-1}),\qquad
    T^t(zI_m)T=z\,{\mathrm{diag}}(G_n,m).\]
  Taking the sum-zero factors, as in \cite[Proposition 2.11]{Yang2025}, yields
  \[\gamma_N =
    \begin{pmatrix}
      x I_n & y G_n^{-1} \\
      z G_n & w I_n
    \end{pmatrix}.
  \]
  Since \(G_n^{-1}=I_n-J_n/m\) and \(m^2\mid y\), the block \(b_N=yG_n^{-1}\) is an integral homomorphism. It is an isogeny because \(b_NG_n=yI_n\) is an isogeny.

  By \Cref{thm:equivalence}, there is a simple semihomogeneous vector bundle \(P_0\) realizing \(\gamma_N\). Orlov's \cite[Corollary 3.4]{Orlov2002} gives
  \[\mathcal Q\cong(t_u\times\id)_*P_0\otimes{\mathrm{pr}}_{N_M}^*L[\ell]
    \quad\text{for }u\in N_M,\ L\in\Pic^0(N_M),\ \ell\in\Z.\]
  Set \(P=\mathcal Q[-\ell]\). It is a locally free simple semihomogeneous equivalence kernel. Shifting the equivariance isomorphisms of \(\mathcal Q\) and taking their degree-zero cohomology gives a diagonal \(S_m\)-linearization of \(P\). The cocycle identities are preserved. Thus the induced equivariant functor is also an equivalence.

  Translation and tensoring by \(L\in\Pic^0\) preserve the divisor slope, so \(\mu(P)=\mu(P_0)\). The slope formula in \Cref{thm:equivalence} gives the symmetric rational homomorphism
  \[\begin{pmatrix}
      G_n y^{-1}x&-G_n y^{-1}\\
      -(y^{-1})^{\vee}G_n&G_n wy^{-1}
    \end{pmatrix}.\]
  For \(F\), the corresponding blocks are obtained by replacing \(G_n\) with \(I_m\). Pullback by \(q_M\times q_A\) multiplies each block by \(T^t\) on the left and \(T\) on the right. Since \(T^tT={\mathrm{diag}}(G_n,m)\), restriction to the sum-zero factors gives precisely the displayed matrix. In particular,
  \[\mu(P)=\iota^*(q_M\times q_A)^*\mu(F),\]
  where \(\iota(u,v)=((u,0),(v,0))\) is the inclusion of \(N_M\times N_A\).

  Finally, the eigenvalues of \(G_n=I_n+J_n\) are \(1\), with multiplicity \(n-1\), and \(n+1=m\). Hence \(\det G_n=m\). On \((M^{\vee})^n\), a product of abelian surfaces, this gives \(\deg G_n=(\det G_n)^4=m^4\). Using \(b_NG_n=yI_n\), \(\deg y=r^2\), and the rank formula in \Cref{thm:equivalence}, we obtain
  \[\rk(P)^2=\rk(P_0)^2=\deg b_N
    =\frac{(\deg y)^n}{\deg G_n}=\frac{r^{2n}}{m^4},\qquad
    \rk(P)=\frac{r^n}{m^2}.\]
\end{proof}

Let \(\rho\) denote the linearization of \(P\). The equivariant equivalence \(\Phi_{P,\rho}\) is the lift obtained in \Cref{prop:rk-p}. For the kernel calculation, its \(G\times G\)-equivariant kernel is the inflation \(\operatorname{Inf}_{G_\Delta}^{G\times G}(P,\rho)\) \cite[\S1.3, Lemma 5]{Ploog2007}. Conjugating by the two McKay equivalences in \Cref{prop:kummer-mckay} yields an equivalence with kernel
\begin{equation}
\label{eq:kernel}
E_{\rho} = {\mathrm{pr}}^*_M \cO(\delta_M) \otimes (h_M \times h_A)_*^{G_{\Delta}}((b_M \times b_A)^* P \otimes {\mathrm{sgn}}).
\end{equation}
Here the group is diagonal in \(G \times G\). This formula follows from equivariant convolution. There are no higher pullbacks of \(P\). Since \(P\) is locally free, the remaining pushforward is finite.

Pullback preserves local freeness. Pushforward under a finite locally free map preserves locally free modules, and taking invariants under the group action gives a direct summand. Thus \(E_{\rho}\) is a vector bundle. Over the free locus, the product has degree \((m!)^2\) and diagonal invariants divide rank by \(m!\), so
\[R\coloneq \rk E_{\rho} = m!p.\]

\begin{lemma}
\label{lem:c1}
  Let \(h:T\to Z\) be a finite flat Gorenstein morphism of complex varieties, with \(Z\) smooth, and let \(V\) be a vector bundle of rank \(p\) on \(T\). Then
  \[c_1(h_*V)=h_*c_1(V)-\frac{p}{2}h_*c_1(\omega_h)
    \quad\text{in }H^2(Z,\QQ).\]
\end{lemma}

\begin{proof}
  Put \(D_h=\det(h_*\cO_T)\). For a line bundle \(L\) on \(T\), its \emph{norm} \(\operatorname{Nm}_h(L)\) is the line bundle on \(Z\) obtained by taking the algebraic norms of its transition functions \cite[Tags 0BD2, 0BCY]{Stacks2026}. Equivalently,
  \[\operatorname{Nm}_h(L)\cong\det(h_*L)\otimes D_h^{-1}.\]
  The determinant--norm identity \cite[\S7.1, equation (7.1.1)]{Deligne1987}, applied to \(V\) and \(\cO_T^{\oplus p}\), gives
  \[\det(h_*V)\cong\operatorname{Nm}_h(\det V)\otimes D_h^{\otimes p}.\]

  The norm--divisor identity \cite[Tag 02RT]{Stacks2026} and the divisor description of the first Chern class \cite[Tag 02SJ]{Stacks2026} yields
  \[c_1(\operatorname{Nm}_h(L))=h_*c_1(L).\]
  This interpretation also applies when \(T\) is singular. Taking first Chern classes in the determinant identity now gives
  \[c_1(h_*V)=h_*c_1(V)+p\,c_1(D_h).\]

  The Gorenstein hypothesis makes \(\omega_h\) a line bundle, and finite duality gives \(h_*\omega_h\cong(h_*\cO_T)^\vee\). Applying the preceding identity with \(V=\omega_h\) and \(p=1\), we obtain
  \[-c_1(D_h)=h_*c_1(\omega_h)+c_1(D_h),\qquad
    c_1(D_h)=-\frac12h_*c_1(\omega_h).\]
  The formula follows from the substitution.
\end{proof}

For \(\lambda \in H^2(A)\), let \(\theta_A(\lambda) = \sum_i {\mathrm{pr}}_i^* \lambda|_{N_A}\). Fix a point \(z \in K(M)\) in the free locus, with a lift \(x \in N_M\). Set \(P_x = P|_{\{x\} \times N_A}\). By \cref{prop:rk-p}, we have
\[\frac{c_1(P_x)}{p} = \frac{c}{r}\theta_A(H).\]
On the correspondence, \(b_A^* \theta_A(H) = h_A^* H\). Hence,
\[\frac{c_1(b_A^* P_x)}{p} = \frac{c}{r} h_A^* H.\]
Then equation \eqref{eq:kernel} identifies \(E_z\) with \(h_{A*}b_A^{*}P_x\), up to tensoring by a constant one-dimensional vector space. Hence, by
\Cref{lem:c1,prop:kummer-mckay}, we have
\[B_A \coloneq \frac{c_1(E_z)}{R} = \frac{c}{r}H - \frac{1}{2} \delta_A.\]

The same computation at a free point \(w \in K(A)\), where \eqref{eq:kernel} gives an additional \(\delta_M\) factor, yields
\[B_M \coloneq \frac{c_1(E_w)}{R} = \frac{k}{r}\hat{H} + \frac{1}{2}\delta_M.\]
Since \(H^1(K(M)) = H^1(K(A)) = 0\), the K\"unneth formula gives
\[H^2(K(M) \times K(A), \QQ) = {\mathrm{pr}}_M^* H^2(K(M), \QQ) \oplus {\mathrm{pr}}_A^* H^2(K(A), \QQ).\]

So
\[\frac{c_1(E_\rho)}{R} = {\mathrm{pr}}^*_M B_M + {\mathrm{pr}}_A^*B_A.\]

The same calculation applies to linearization \(\rho \otimes {\mathrm{sgn}}\).

\section{Deformation and Equivalence}

Let \(\psi_{E_\rho}\) be the positive cone normalized rational isometry extracted from \([\kappa(E_{\rho})\sqrt{{\mathrm{td}}}]_*\).

Define vector bundles on the full product \(N_M\times N_A\) by
\[P_{\sigma}=(\mathrm{id}_{N_M}\times\sigma)^*P,\qquad
  W=\bigoplus_{\sigma\in G}P_{\sigma}.
\]
For \(S=M,A\), let \(N_S^{\circ}\subset N_S\) be the ordered locus of distinct points and let \(K(S)^{\circ}\) be its quotient by \(G\), identified with the corresponding open subset of \(K(S)\). Write
\[\pi_S:N_S^{\circ}\longrightarrow K(S)^{\circ}\]
for the finite \'etale ordered cover. The complement of \(N_S^{\circ}\) has codimension \(2\): for \(m\geq3\), each difference homomorphism \(N_S\to S\), \((x_i)\mapsto x_i-x_j\) with \(i\ne j\), is surjective.

The free-locus description of \eqref{eq:kernel} gives an isomorphism of ordinary vector bundles
\begin{equation}
  \label{eq:ordered-kernel}
  (\pi_M\times\pi_A)^*\bigl(E_\rho|_{K(M)^{\circ}\times K(A)^{\circ}}\bigr)
  \cong W|_{N_M^{\circ}\times N_A^{\circ}}.
\end{equation}
Thus \(W\) is defined globally on \(N_M\times N_A\), and its comparison with \(E_\rho\) is made after restriction and pullback to the ordered free loci.

Write \(c_1(P_\sigma)/p = a_M + a_{A} + \chi_{\sigma}\). The first two factors are independent of \(\sigma\), and \(\sum_{\sigma} \chi_\sigma = 0\) by permutation. Each \(P_{\sigma}\) is projectively flat, so
\[\frac{\kappa(W)}{R} = \frac{1}{m!} \sum_{\sigma} e^{\chi_{\sigma}}.\]

\begin{lemma}
\label{lem:kappa2}
With \(\chi\) in \eqref{eq:chi},
    \[\left[\frac{\kappa(W)}{R}\right]_2 = \frac{1}{2m!}\sum_{\sigma \in G} \chi_{\sigma}^2 = \frac{1}{n}(\theta_M \otimes \theta_A)(\chi^2/2).\]
Here \(\theta_S(\lambda)=\sum_{k=1}^m{\mathrm{pr}}_{S,k}^*\lambda\) for \(S=M,A\), where \({\mathrm{pr}}_{S,k}:N_S\to S\) is the \(k\)-th coordinate projection.
\end{lemma}

\begin{proof}
  For \(S=M,A\) and \(\eta\in H^1(S,\R)\), the sum-zero condition gives \(\sum_k{\mathrm{pr}}_{S,k}^*\eta=0\). Thus
  \[H^1(N_S,\R)\cong H^1(S,\R)\otimes V,\qquad
    V=\{(z_k)\in\R^m:\textstyle\sum_kz_k=0\},
  \]
  identifying \(V\) with its dual by the standard inner product. Here \(\eta\otimes v\) corresponds to \(\sum_k v_k{\mathrm{pr}}_{S,k}^*\eta\). Coordinate permutations act on \(V\) through the standard representation, of dimension \(n=m-1\).

  Write \(\chi=\sum_{a,b}J_{ab}u_a\cup w_b\). For an orthonormal basis \(e_i=(e_{ik})\) of \(V\), set \(u_{a,i}=u_a\otimes e_i\) and \(w_{b,i}=w_b\otimes e_i\), viewing these classes on the product. The projection identity \(\sum_i e_{ik}e_{i\ell}=\delta_{k\ell}-1/m\), together with the vanishing coordinate sums, gives
  \[\sum_i u_{a,i}\cup u_{c,i}=\theta_M(u_a\cup u_c),\qquad
    \sum_j w_{b,j}\cup w_{d,j}=\theta_A(w_b\cup w_d).
  \]
  Similarly, \Cref{prop:rk-p} identifies the mixed slope of \(P\) as
  \[\chi_P=\sum_{a,b}J_{ab}\sum_k
    {\mathrm{pr}}_{M,k}^*u_a\cup{\mathrm{pr}}_{A,k}^*w_b
    =\sum_{a,b,i}J_{ab}u_{a,i}\cup w_{b,i}.
  \]
  Let \(R(\sigma)\) be the matrix of \((\sigma^{-1})^*\) on \(V\). Pullback by \(\mathrm{id}\times\sigma\) permutes the target factor, so
  \[\chi_\sigma=\sum_{a,b,i,j}J_{ab}R(\sigma^{-1})_{ji}
    u_{a,i}\cup w_{b,j}.
  \]
  The standard representation is orthogonal and irreducible, so Schur orthogonality gives
  \begin{equation}
    \label{eq:avg}
    \frac1{m!}\sum_{\sigma\in G}R(\sigma)_{ji}R(\sigma)_{\ell k}
    =\frac1n\delta_{ik}\delta_{j\ell}.
  \end{equation}
  Reindexing by \(\sigma^{-1}\) and applying \eqref{eq:avg} to the square gives
  \[
    \begin{aligned}
      \frac1{m!}\sum_{\sigma\in G}\chi_\sigma^2
      &=-\frac1n\sum_{a,b,c,d}J_{ab}J_{cd}
        \theta_M(u_a\cup u_c)\cup\theta_A(w_b\cup w_d)\\
      &=\frac1n(\theta_M\otimes\theta_A)(\chi^2).
    \end{aligned}
  \]
  The minus sign comes from interchanging the degree-one classes \(w_{b,j}\) and \(u_{c,i}\). The same reordering in \(\chi^2\) explains the second equality. Taking the degree-four term \(\chi_\sigma^2/2\) of \(e^{\chi_\sigma}\) proves the formula for \(\kappa_2(W)/R\).
\end{proof}

We use \Cref{lem:kappa2} to determine \(\psi_{E_\rho}\). Recall the orthogonal decompositions
\[H^2(K(S),\QQ)=H^2(S,\QQ)\oplus\QQ\delta_S,\qquad S=M,A.\]
We first determine the \(H^{2}(A,\QQ)\)-component of the image of a class in \(H^{2}(M,\QQ)\), and then use the isometry condition to find the remaining components.

By \Cref{lem:degree-reversing}, \(\psi_{E_\rho}\) is a nonzero scalar multiple of
\[x\longmapsto[\kappa(E_\rho)\sqrt{\td}]_*\bigl(x\,c_2(K(M))^{n-1}\bigr).\]
The input has degree \(4n-2\), and pushforward integrates over \(K(M)\), of dimension \(2n\). Thus only the component of the correspondence in
\[H^2(K(M),\QQ)\otimes H^2(K(A),\QQ)\]
contributes to the degree-two output. Since \(\kappa_1(E_\rho)=0\) and \(\sqrt{\td}\) has no component in this tensor product, this component comes entirely from \(\kappa_2(E_\rho)\).

The action of this mixed class is computed by pairing its source factor with \(x\). The Fujiki relation gives
\[\int_{K(M)}x\,u\,c_2(K(M))^{n-1}=C_n(x,u)_{K(M)},\qquad C_n\neq0,\]
for \(x,u\in H^2(K(M),\QQ)\), as in the proof of \cite[Lemma 3.7]{Markman2024}. On the surface summand, the BBF form agrees with the intersection pairing on \(M\). Moreover, a surface class \(x\) is orthogonal to \(\delta_M\), so a term whose source factor is \(\delta_M\) contributes nothing.

To recover the mixed class from \Cref{lem:kappa2}, we use \eqref{eq:ordered-kernel}. For \(S=M,A\), restriction to the free locus followed by pullback to the ordered cover has kernel
\[\ker\bigl(H^2(K(S),\QQ)\longrightarrow H^2(N_S^{\circ},\QQ)\bigr)=\QQ\delta_S.\]
Hence the mixed class is determined by its pullback up to terms containing \(\delta_M\) or \(\delta_A\). For a surface input, the terms containing \(\delta_M\) vanish by orthogonality, while those containing \(\delta_A\) affect only the exceptional component of the output. Therefore \Cref{lem:kappa2} determines the component in \(H^2(A,\QQ)\): it is a nonzero scalar multiple of \([r\chi^2/2]_*|_{H^2(M)}\), and hence of \(\psi_{\cU}\) from \eqref{eq:surface-isometry}.

It follows that there are a nonzero scalar \(a\in\QQ\) and a linear form \(\ell:H^2(M,\QQ)\to\QQ\) such that
\[\psi_{E_\rho}(x)=a\,\psi_{\cU}(x)+\ell(x)\delta_A\qquad\text{for }x\in H^2(M,\QQ).\]
Both \(\psi_{E_\rho}\) and \(\psi_{\cU}\) are isometries. Using the orthogonality of the target summands and \(\delta_A^2=-2m\), we obtain
\[(a^2-1)(x,y)_M=2m\,\ell(x)\ell(y)\qquad\text{for }x,y\in H^2(M,\QQ).\]
If \(a^2\neq1\), the left-hand side is a bilinear form of rank six, whereas the right-hand side has rank at most one. Thus \(a^2=1\), and the same identity then forces \(\ell=0\). Since \(\psi_{\cU}(\hat H)=H\) and \(\hat H,H\) lie in the chosen positive components, the positive-cone normalization of \(\psi_{E_\rho}\) selects \(a=1\).

Finally, an isometry mapping the surface summands onto each other also maps their orthogonal complements onto each other. Since \(\delta_M^2=\delta_A^2=-2m\), the image of \(\delta_M\) is \(\epsilon\delta_A\) for some \(\epsilon\in\{1,-1\}\). We have therefore proved
\begin{equation}
  \label{eq:e-rho}
  \psi_{E_\rho} = \psi_{\cU} \oplus (\epsilon \id_{\delta}), \quad \epsilon \in \{1, -1\}.
\end{equation}

For either Kummer variety, we write
\[T_A = {\mathrm{BKR}}_A \circ (- \otimes {\mathrm{sgn}}) \circ {\mathrm{BKR}}_A^{-1}.\]
This is an involution and fixes the skyscraper sheaves in the free locus. Since \((h_{A*}\cO_{\Gamma_A})^G = \cO_{K(A)}\), the trivial equivariant structure sheaf maps to \(\cO_{K(A)}\). Finite duality gives
\[(h_{A*}\omega_{h_A})^G \cong ((h_{A*}\cO_{\Gamma_A})^G)^{\vee} \cong \cO_{K(A)}.\]
Using \(\omega_{h_A} \cong {\mathrm{sgn}} \otimes h_A^*\cO(\delta_A)\) from \Cref{prop:kummer-mckay} and the projection formula, we have
\[T_A(\cO) = \cO(-\delta_A).\]
Let \(\delta = \delta_{A}\). We write \(\rho_{\delta}\) for the reflection fixing \(\delta^{\perp}\) and sending \(\delta\) to \(-\delta\), that is,
\[\rho_{\delta}(x) =  x - \frac{2(x, \delta)}{(\delta, \delta)}\delta.\]
On \(\tilde{H}\), write \(\mathsf{e}_{\lambda}\) for the isometry
\[\mathsf{e}_{\lambda}(\alpha) = \alpha + \lambda + \frac{1}{2}q(\lambda)\beta, \mathsf{e}_{\lambda}(x) = x + (x, \lambda) \beta, \mathsf{e}_{\lambda}(\beta) = \beta.\]
Under the LLV representation, this isometry induces the exponential cup product operator \(\mathsf{e}_{\lambda}(u)=e^{\lambda}\cup u\) defined above.

\begin{lemma}
  \label{lem:ta}
    The action of \(T_A\) on the rank and first Chern class is
    \[(r, \lambda) \to (r, \rho_{\delta}(\lambda) - r \delta),\]
    and on the extended orthogonal group its projective class is \([\mathsf{e}_{-\delta}\rho_{\delta}]\)
\end{lemma}

\begin{proof}
  The extended LLV space is \(\tilde{H}\). The subrepresentation \(SH^{*}\) generated by \(H^{2}\) is the irreducible trace-free symmetric power representation of degree \(n\). Under its standard identification, \(1\) corresponds to \(\alpha^{n}/n!\) and the point class is a nonzero multiple of \(\beta^{n}\). By \Cref{thm:LLV}, the cohomological action of \(T_A\) preserves \(SH^{*}\) and is a Mukai isometry.

  The Mukai pairing restricted to \(SH^{*}\) is nondegenerate, so this action also preserves \(SH^{*\perp}\). Since the grading operator belongs to the LLV algebra, this complement is graded. Its degree-zero and degree-two spaces vanish because \(SH^{*}\) contains all of \(H^0\) and \(H^2\). Therefore \(SH^{*\perp}\) does not affect the rank and first Chern class.

  As \(T_{A}\) fixes the point class, its extended isometry preserves \(\QQ \beta\). An isometry preserving \(\QQ \beta\) has the form
  \[\beta \to t \beta, x \to A_{0}x + l_{0}(x)\beta, \alpha \to t^{-1}\alpha + v + s_{0}\beta.\]
  It cannot decrease the weight, so its symmetric power preserves the descending cohomological filtration. Therefore the higher degree part of \(SH^{*}\) also does not affect the rank and first Chern class.

  The rank is preserved, since
  \[\rk T(F) = \chi(T(F), \cO_{x}) = \chi(F, \cO_{x}) = \rk F.\]
  So by linearity, \(T_{A}\) maps
  \[(u, \lambda) \to (u, B\lambda + ub),\]
  for some linear map \(B: H^{2}(K(A)) \to H^{2}(K(A))\) and \(b \in H^{2}(K(A))\). Now using \(T(\cO) = \cO(-\delta)\), we obtain \(b = - \delta\). Now it remains to show that \(B\) is the exceptional reflection. First, notice that \(B\) is a similitude,
  \[(Bx, By) = \mu(x, y).\]
  On \(K(A)^{\circ}\), \(T_A\) is tensoring by the line bundle \(L_{\mathrm{sgn}}\) associated to the sign character of the ordered \(S_m\)-cover. Since \(L_{\mathrm{sgn}}^{\otimes 2} \cong \cO\), its first Chern class vanishes in rational cohomology. So \(T_A\) acts trivially on rational cohomology of the free locus. Since
  \[\ker (H^{2}(K(A), \QQ) \to H^{2}(K(A)^{\circ}, \QQ)) = \QQ \delta,\]
  we have \(B\lambda - \lambda \in \QQ \delta\). Take \(\sigma \in H^{2,0}(K(A))\). The action preserves Hodge types, so \(B \sigma - \sigma\) has type \((2, 0)\). So by the above formula, it must be zero, that is, \(B \sigma = \sigma\). The same holds for \(\bar{\sigma} \in H^{0,2}(K(A))\). Since \((\sigma, \bar{\sigma}) \neq 0\), we have
  \[(\sigma, \bar{\sigma}) = (B\sigma, B\bar{\sigma})\]
  That is, \(\mu = 1\), and \(B\) is an isometry. So \(B\delta = a \delta\) with \(a = \pm 1\). If \(x \perp \delta\), write \(Bx = x + b \delta\). Then
  \[0 = (Bx, B \delta) = ab(\delta, \delta).\]
  So \(b = 0\), and \(B\) fixes \(\delta^{\perp}\). Applying \(T\) twice to the pair \((1, 0)\), we obtain \((1, -B\delta - \delta)\). Since \(T^{2} = \id\), \(B\) negates \(\delta\). Therefore \(B = \rho_{\delta}\).

  Let \(g: \tilde{H}(K(A)) \to \tilde{H}(K(A))\) be a representative of the projective extended action and \(T|_{SH^{*}} = a \Sym_{0}^{n}(g)\). Since it preserves the point class, we have \(g(\beta) = t\beta\). \(g\) acts on the point class by multiplying \(at^{n}\). On the other hand, the degree zero component is multiplied by \(at^{-n}\). Since \(T\) fixes the point class and preserves the rank, we have \(t^{2n} = 1\). That is, \(t = \pm 1\). Negating if necessary, we take \(t = 1\) and hence \(a = 1\). Then
  \[\beta \to \beta, x \to A_{0}x + (A_{0}x, v) \beta, \alpha \to \alpha + v + \frac{1}{2}q(v)\beta.\]
  Since \(T(1) \equiv 1 + v \pmod{H^{\geq 4}}\) and \(T(x) = A_{0}x \pmod{H^{\geq 4}}\), the induced low-degree action is
  \[(u, \lambda) \to (u, A_{0}\lambda + uv).\]
  Comparing with the previous calculation of the first Chern class gives \(A_{0} = \rho_{\delta}\) and \(v = - \delta\), hence \(g = \mathsf{e}_{-\delta}\rho_{\delta}\). This finishes the proof.
\end{proof}

Set \(E=E_\rho\) and \(E'=E_{\rho\otimes{\mathrm{sgn}}}\). Changing the linearization gives \(\Phi_{E'}=T_A\Phi_E\). The kernels \(E\) and \(E'\) have the same source and target determinant slopes \(B_M,B_A\), as computed above. Since \(\kappa(E)=\ch(E)e^{-{\mathrm{pr}}_M^*B_M-{\mathrm{pr}}_A^*B_A}\), the projection formula gives
\[\phi_E=\mathsf{e}_{-B_A}\Phi_E^H\mathsf{e}_{-B_M}.\]
Here \(\phi_E\) is the normalized cohomological transform defined in \Cref{lem:degree-reversing}, with \(X=K(M)\) and \(Y=K(A)\), and \(\Phi_E^H\) is the ordinary cohomological transform, and \(\mathsf{e}_a\) acts on cohomology by multiplication by \(e^a\). Applying the same normalization to \(E'\), we obtain
\[
  \begin{aligned}
    \phi_{E'}
    &=\mathsf{e}_{-B_A}T_A^H\Phi_E^H\mathsf{e}_{-B_M}\\
    &=(\mathsf{e}_{-B_A}T_A^H\mathsf{e}_{B_A})\phi_E.
  \end{aligned}
\]
By \Cref{lem:ta}, we have
\[
  \begin{aligned}
    \mathsf{e}_{-B_A}(\mathsf{e}_{-\delta_A}\rho_{\delta_A})\mathsf{e}_{B_A}
    &=\mathsf{e}_{-B_A-\delta_A+\rho_{\delta_A}(B_A)}\rho_{\delta_A}\\
    &=\rho_{\delta_A}.
  \end{aligned}
\]
Therefore the positive-cone normalized isometries satisfy
\[\psi_{E_{\rho\otimes{\mathrm{sgn}}}}=\rho_{\delta_A}\circ\psi_{E_\rho},\]
so the two kernels have opposite exceptional signs in \eqref{eq:e-rho}, with the same surface isometry \(\psi_{\cU}\). We may therefore choose one \(E\) with
\[\psi_E = \psi_{\cU} \oplus \id_{\delta}.\]
For sufficiently small positive rational \(\epsilon\), the ample class \(\hat{H} - \epsilon \delta_M\) and \(H - \epsilon \delta_A\) are matched by \(\psi_E\).

Let \(Z = K(M) \times K(A)\) and \(N = \dim_{\C}Z = 4n\). Choose an integral multiple \(L\) of the pullback of the quotient product polarization \(\hat{H} + H\). It is big and nef, and a Hilbert-Chow exceptional divisor has \(L\)-degree zero.

On \(N_M \times N_A\) the bundles \(P_{\sigma}\) are stable for every ample polarization. And for different \(\sigma\), the numerical slopes agree. Also, they are pairwise non-isomorphic: an isomorphism of two kernels, composed with the inverse of \(\Phi_P\), would identify the two permutation functors. Applying them to skyscraper sheaves identifies the permutations on \(N_A\). The action of \(S_m\) on \(N_A\) is faithful for \(m \geq 3\). Thus \(W = \oplus_{\sigma}P_{\sigma}\) is a direct sum of distinct stable bundles with the same numerical slope.

\begin{lemma}
  \label{lem:E-stable}
  The bundle \(E\) is slope-stable with respect to \(L\).
\end{lemma}

\begin{proof}
  Let \(0 \neq F \subset E\) be saturated with \(\rk F < R\). Restrict it to \(K(M)^{\circ}\times K(A)^{\circ}\) and pull it back by \(\pi_M\times\pi_A\). Via \eqref{eq:ordered-kernel}, this is a subsheaf of \(W|_{N_M^{\circ}\times N_A^{\circ}}\). Extend it and take its saturation in the globally defined bundle \(W\) on \(N_M\times N_A\), obtaining a coherent saturated subsheaf \(\tilde{F}\subset W\). The first Chern class on the smooth torus does not change on deleting a codimension \(2\) subset, so on \(Z\) the ambiguity in extending a class from the free locus is supported on the exceptional divisors, whose intersection number with \(L^{N - 1}\) is zero. Therefore, the projection formula for finite quotient identifies the slope difference upstairs with the downstairs difference multiplied by the positive covering degree \((m!)^{2}\). Now polystability of \(W\) gives \(\mu_{L}(F) \leq \mu_{L}(E)\). If equality held, a saturated equal-slope subsheaf \(\tilde{F}\) would be a direct summand. The \(G \times G\) action on the labels of \(P_{\sigma}\) is transitive, and \(\tilde{F}\) is invariant under \(G \times G\) action. So \(\tilde{F}\) must be \(W\), which contradicts our choice of \(\rk F < R\). Hence the inequality is strict.
\end{proof}

\begin{lemma}
  \label{lem:eta-stable}
    If a vector bundle is stable with respect to an integral big and nef class \(L\), it is stable with respect to \(L + \eta L'\) for every sufficiently small \(\eta > 0\), for a fixed integral ample divisor \(L'\).
\end{lemma}

\begin{proof}
  Choose \(k_{0}\) such that \(E^{\vee}(k_{0}L')\) is globally generated. Dualizing its evaluation map embeds \(E\) into a sum of copies of \(\cO(k_{0}L')\). For a saturated subsheaf \(F\) of rank \(r' < R\), taking the \(r'\)th exterior power and extending over codimension \(2\) gives \(H^{0}(Z, \cO_{Z}(k_{0}r'L') \otimes (\det F)^{-1}) \neq 0\), and therefore
  \[D_{F} \coloneq k_{0}r'L'  - c_{1}(F)\]
  is effective. For every \(j\), effectiveness and nefness give a uniform upper bound on
  \[a_{j}(F) = \left( \frac{c_{1}(F)}{r'} - \frac{c_{1}(E)}{R} \right) L^{N - 1 - j}L^{\prime j} = (k_{0}L' - \frac{c_{1}(E)}{R}) \cdot L^{N - 1 - j} L^{\prime j} - \frac{1}{r'} D_{F} \cdot L^{N - 1 - j} L^{\prime j},\]
  independent of \(F\). Moreover, for \(j = 0\) the value is a strictly negative rational number at most \(-1/(R(R - 1))\). Hence in
  \[\mu_{L + \eta L'}(F) - \mu_{L + \eta L'}(E) = \sum_{j = 0}^{N - 1} \binom{N - 1}{j} \eta^{j}a_{j}(F)\]
  the first term is \(\leq -1/(R(R - 1))\), and all other coefficients have fixed uniform upper bounds. So a sufficiently small positive \(\eta\) makes the sum of their positive contributions smaller than \(1/(R(R - 1))\), simultaneously for every \(F\).
\end{proof}

We choose a compatible ample class \(L'\), for instance an integral multiple of \((\hat{H} - \epsilon \delta_{M}) + (H - \epsilon \delta_{A})\). The class \(L + \eta L'\) is compatible as well. By \Cref{lem:eta-stable}, \(E\) is stable with respect to \(\Omega_{0} \coloneq L + \eta_{0} L'\) for sufficiently small \(\eta_{0}>0\). Openness of slope stability gives a nonempty open subcone of matched K\"ahler classes on which \(E\) is stable \cite[Corollary 3.5]{Delloque2025}. To apply \Cref{prop:azumaya-continue}, it remains to check that the transported \(c_{2}(\cEnd E)\) is of type \((2,2)\) along the compatible deformations.

Fix markings of \(K(M)\) and \(K(A)\), let \(\psi\) be the correspondence induced by \(\psi_E\) and \(\cM_\psi^0\) be the component containing this pair. Write \(X=K(M)\), \(Y=K(A)\), and set \(v_E=\kappa(E)\sqrt{\td_{X\times Y}}\), \(\phi_E=[v_E]_*\). For our locally free equivalence kernel \(E\), \Cref{lem:degree-reversing} gives LLV equivariance and
\[\Ad_{\phi_E}|_{\bar{\mathfrak g}_X}=\Ad_{\psi_E}.\]
The argument of \cite[proof of Proposition 5.22, Step 1]{Markman2024} therefore shows that, along every path of compatible pairs \((X_t,Y_t)\) in \(\cM_\psi^0\), the parallel transport \(v_t\) of \(v_E\) remains a Hodge class:
\[v_t\in\bigoplus_{k\geq0}H^{k,k}(X_t\times Y_t).\]
That argument uses only LLV equivariance and Hodge compatibility, and does not assume the vanishing of odd cohomology, so it applies to our generalized Kummer type factors. We refer to the cited proof for the details.

Since the Todd classes remain Hodge classes, the transported \(\kappa_2(E)\) has type \((2,2)\). The identity
\[c_2(\cEnd E)=2R\,c_2(E)-(R-1)c_1(E)^2=-2R\,\kappa_2(E)\]
then shows that the transported \(c_2(\cEnd E)\) remains of type \((2,2)\), as required by \Cref{prop:azumaya-continue}.

Therefore, we can apply \Cref{prop:azumaya-continue}. Each segment is associated to an actual product hyperk\"ahler metric, and we obtain an Azumaya algebra on each pair. We write \(X_{1}\) and \(X_{2}\) for the source and target in the deformation, and \(F\) for a locally free, possibly twisted bundle of rank \(R\) on \(X_{1} \times X_{2}\).

\begin{proposition}[{\cite[Theorems 2.17--2.18]{HS2026}, \cite[Theorem 2.3]{KapustkaKapustka2025}}]
  \label{prop:ext-continue}
  Let \(q = (X_{1}, \eta_{1}, X_{2}, \eta_{2}) \in \cM_{\psi}^{0}\) be a projective pair, and let \(F\) be obtained by continuing \(E\) along a generic diagonal twistor path as in \Cref{prop:azumaya-continue}. Then \(F\) induces an exact equivalence
  \[\Phi_F: D^{b}(X_{1}, \beta_{1}) \xrightarrow{\sim} D^{b}(X_{2}, \beta_{2}),\]
  where \(\beta_i \in \Br(X_i)\), and the twisting class of \(F\) is \(-{\mathrm{pr}}_{1}^{*}\beta_{1} + {\mathrm{pr}}_{2}^{*}\beta_{2}\).
\end{proposition}

\begin{proof}
  The bundle \(E\) is untwisted and induces a Fourier-Mukai equivalence. By \Cref{thm:LLV}, its cohomological transform is LLV-equivariant. Its initial Brauer classes are zero, and stability on an open subcone of matched K\"ahler classes was proved above. So the hypotheses of \cite[Theorem 2.17]{HS2026} hold for the continuation in \Cref{prop:azumaya-continue}. Now \cite[Theorem 2.18]{HS2026} gives the equivalence.
\end{proof}

Thus every projective pair in \(\cM_{\psi}^{0}\) admits a twisted derived equivalence. We compute the source and target Brauer classes from the transported determinant slopes.

\section{Wall crossing}
\label{sec:wall-crossing}
We start from a vector bundle \(E_{\rho}\) on \(K(M) \times K(A)\) and deform \(K(M) \times K(A)\) to \(Y \times X\). The topological twisting class of the projective bundle is \(-c_1(E)\) modulo \(R\), and it transports with the smooth projective bundle along the whole twistor path. By \cite[Propositions 4.2--4.3]{Caldararu2002}, its image in the holomorphic Brauer group at the endpoint is \(\left[-\tau(c_1(E))/R\right]_{Y \times X}\), where \(\tau\) is parallel transport. Changing a smooth lift by a line bundle changes this rational class by an integral class, whose Brauer class is zero. So the same formula holds through the twistor joints, even when the intermediate Brauer classes are nonzero. Now let us apply this to the product. Let \(B_{Y} = \tau_{Y}(B_{M})\) and \(B_{X} = \tau_{X}(B_{A})\) the transported classes on \(Y\) and \(X\). We have
\[\gamma = -p_{Y}^{*}[B_{Y}]_{Y} - p_{X}^{*}[B_{X}]_{X}.\]
Let \(F\) denote the endpoint twisted kernel, and set
\[\beta = [B_{Y}]_{Y}, \alpha = -[B_{X}]_{X}.\]
Then the twist of \(F\) is \(-p_{Y}^{*}\beta + p_{X}^{*}\alpha\), so \(F\) induces a morphism
\[\Phi_{F}: D^{b}(Y, [B_{Y}]_{Y}) \to D^{b}(X, -[B_{X}]_{X}).\]
By \Cref{prop:ext-continue}, it is an equivalence. If \([B_{X}]_{X}\) vanishes, its target is \(D^{b}(X)\). If two paths end at the same marked source \(Y\), the transported degree \(2\) class \(B_{M}\) is the same, and so is its Brauer class.

\begin{lemma}
\label{lem:lattice-1}
    Let \(S \subset \Lambda = U^{\oplus 3} \oplus \Z \delta\) be a primitive sublattice of rank at least \(2\), and let \(\lambda: S \to \Z\) be a nonzero linear form. There exists \(D \in S\) with \({\mathrm{div}}_{\Lambda}(D) = 1\) and \(\lambda(D) \neq 0\).
\end{lemma}

\begin{proof}
    Choose two members \(u_1\), \(u_2\) of a basis of \(S\), with \(\lambda(u_1) \neq 0\). Because \(S\) is primitive, the reduction \(S/pS \to \Lambda/p\Lambda\) is injective for every prime \(p\). For \(p | 2m\), the radical of the bilinear form on \(\Lambda/p\Lambda\) is exactly the one-dimensional line generated by \(\delta\).

    At most one residue \(z \in \mathbb F_p\) can satisfy \(u_1 + zu_2 \in {\mathrm{rad}}(\Lambda/p\Lambda)\). For each prime dividing \(2m\), choose a residue \(z_p\) avoiding this value. By the Chinese remainder theorem there are infinitely many \(z\) with residue \(z_p\) for such primes.

    The vector \(D_0 = u_1 + zu_2\) is primitive in \(S\), hence primitive in \(\Lambda\). Choose \(z\) such that \(\lambda(u_1) + z \lambda(u_2) \neq 0\). So no prime \(p \mid 2m\) divides \({\mathrm{div}}_{\Lambda}(D_0)\) by construction. If \(p \not \mid 2m\), the lattice is nondegenerate modulo \(p\), and the primitive vector \(D_0\) is nonzero modulo \(p\), so some pairing with it is nonzero modulo \(p\). Such a prime cannot divide its divisibility either. Thus the divisibility is \(1\).
\end{proof}

\begin{lemma}
\label{lem:lattice-2}
    Let \(X\) be a projective smooth hyperk\"ahler variety of generalized Kummer type, let \(w \in \NS(X)\) be a primitive negative wall class, and let \(K > 0\). There are an integer \(C > 0\) and classes \(D \in \NS(X)\) in the desired positive component with
    \[{\mathrm{div}}_{\Lambda}(D) = 1, \quad (D, w) = C, \quad q(D) > CK.\]
\end{lemma}

\begin{proof}
    The N\'eron-Severi lattice is primitive, non-degenerate, and has signature \((1, \rho - 1)\). Its rank is at least two because it contains a positive class and the negative class \(w\). So we apply \Cref{lem:lattice-1} with \(\lambda = (-, w)\), negating if necessary, and obtain a \(D_0\) with \(C = (D_0, w) > 0\).

    The rational space \(\NS(X)_{\QQ} \cap w^{\perp}\) contains a positive integral vector \(V\) in the chosen positive component. Let \(D_i = D_0 + iV\). Then \((D_i, w) = C\) and \(q(D_i)\) tends to \(\infty\). Choose \(u \in \Lambda\) such that \((D_0, u) = 1\). Since \(V\) is not proportional to \(D_0\), there is \(v \in \Lambda\), with \((D_0, v) = 0\) and \(b = (V, v) \neq 0\). Write \(a = (V, u)\). If \(n\) is a positive multiple of \(\abs{b}\) and \(e = {\mathrm{div}}(D_n)\), then
    \[e \mid 1 + na,\quad e \mid nb\]
    The first equation gives \({\mathrm{gcd}}(e,n) = 1\) and the second gives \(e \mid b\), while \(b \mid n\) by our choice, so \(e = 1\). For all sufficiently large \(n\), \(D_n\) lies in the required positive component. So taking \(D = D_{k\abs{b}}\) for sufficiently large \(k > 0\) satisfies all the requirements.
\end{proof}

Now by \cref{thm:bound-q}, there is a constant \(B_{m} > 0\) such that every primitive integral MBM class \(v\) of this deformation type satisfies
\[0 < - q(v) \leq B_m.\]
Now choose \(D\) from \Cref{lem:lattice-2} with \(K = 2\sqrt{B_m}\), and set \(g = q(D)/2\). Since \(\Lambda\) is even, \(g\) is a positive integer. Then
\[q(D) = 2g, \quad g > C \sqrt{B_m}.\]
With this \(g\) in \Cref{sec:construct}, the vector
\[h_0 = H - \frac{cg}{2}\delta_A = H - 2m^2g\delta_A\]
has
\[q(h_0) = 2gs - \frac{mc^2g^2}{2} = 2g, \quad {\mathrm{div}}(h_0) = 1.\]
The latter follows from primitivity of \(H\). The class lies in the positive component containing \(H\), since \((h_0, H) = H^2 > 0\). To compare different manifolds, we choose a parallel transport identifying their second integral cohomology. The image of \(h_0\) and \(D\) are both primitive, have square \(2g\), and have trivial discriminant class. So applying \cite[Proposition 3.3]{GHS2009} to \(\Lambda \supset U^{\oplus 2}\), we obtain an element of \(\widetilde{\SO}^{+}(\Lambda)\) carrying one to the other. The monodromy inclusion in \Cref{thm:kummer-monodromy} makes this a monodromy operator. Composing the two parallel transports gives a parallel-transport isometry
\[\rho: H^2(K(A), \Z) \xrightarrow{\sim} H^2(X, \Z), \quad \rho(h_0) = D.\]
Here \(\rho\) is the parallel-transport isometry, whereas \(E_\rho\) is indexed by the linearization of \(P\). We have
\[B_A = \frac{h_0}{cg},\quad \rho(B_A) = \frac{D}{cg} \in \NS(X)_{\QQ}.\]
Thus the Brauer twist class associated to \(\rho(h_0)\) vanishes on \(X\), and also on any birational model \(X'\) with its identified divisor lattice.

\begin{lemma}
\label{lem:divisibility}
    There is an integer \(l\) with \(k =cl\) and \({\mathrm{gcd}}(l,r) = 1\). For every \(\lambda \in H^2(A, \Z)\) and integer \(a\),
    \[a \psi_{\cU}^{-1}(\lambda) \in H^2(M, \Z) \iff g \mid a(H, \lambda).\]
\end{lemma}

\begin{proof}
    The inverse surface equivalence has kernel \(\cU^{\vee}[2]\). Applied to a point of \(A\), its Mukai vector is
    \[(r, -k \hat{H}, \frac{k^2s}{c^2}).\]
    It is integral and primitive. Since \({\mathrm{gcd}}(c, s) = 1\), we must have \(c^2 \mid k^2\) and therefore \(c \mid k\), so we can write \(k = cl\). If a prime divides both \(l\) and \(r\), it divides \((r, -cl\hat{H}, l^2s)\), contradicting primitivity. Hence \(\gcd(l, g) = 1\).

    Use \(f_0^{-1} = \mathsf{e}_{a_M}f^{-1}\mathsf{e}_{a_A}\). For \(\lambda \in H^2(A, \Z)\),
    \[\mathsf{e}_{a_A}(\lambda) = \lambda + \frac{c}{r}(H, \lambda)[{\mathrm{pt}}].\]
    The rank of \(f^{-1}(\lambda)\) is \(-c(H, \lambda)\) and the added point-vector has opposite rank, so the total rank is \(0\). Its degree \(2\)-part is
    \[[f^{-1}(\lambda)]_2 - \frac{ck}{r}(H, \lambda) \hat{H}.\]
    The final multiplication by \(\mathsf{e}_{a_M}\) does not change the degree \(2\)-part, since the rank is \(0\). Therefore
    \[\psi_{\cU}^{-1}(\lambda) = \varepsilon_{\cU}([f^{-1}(\lambda)]_2 - \frac{l}{g}(H, \lambda) \hat{H}).\]

    The first term is integral, and \(\hat{H}\) is primitive. An integral multiple \(a\) of \(\psi_{\cU}^{-1}(\lambda)\) is integral when \(g \mid al(H, \lambda)\), equivalently \(g \mid a(H, \lambda)\).
\end{proof}

Now set \(\Psi = \rho \psi_E\). Write
\[\rho^{-1}(w) = \lambda + t \delta_A, \quad \lambda \in H^2(A, \Z), t \in \Z.\]
Then we have
\[C = (D, w) = (H, \lambda) + 4m^3 gt,\]
and therefore
\[(H, \lambda) \equiv C \pmod g.\]
\begin{proposition}
  \label{prop:no-wall}
    The rational class \(\Psi^{-1}(w)\) is not proportional to a primitive integral MBM class of the deformation type \(K(A)\).
\end{proposition}

\begin{proof}
    Suppose \(\Psi^{-1}(w) = (b/a)w'\) with \(w'\) primitive integral MBM, \(a > 0\), \(b \neq 0\) and \(\gcd(a, b) = 1\). Integrality of \(a\Psi^{-1}(w)\) and \Cref{lem:divisibility} implies \(g \mid aC\). On the other hand, the isometry gives
    \[a^2 q(w) = b^2q(w').\]
    Since \(a^2\) and \(b^2\) are coprime, \(a^2 \mid q(w')\). The norm bound then yields \(a \leq \sqrt{B_m}\). So
    \[0 < aC \leq C \sqrt{B_m} < g,\]
    by our choice of \(g\). This contradicts \(g \mid aC\).
\end{proof}

\begin{theorem}
  Let \(X\) and \(X'\) be birational smooth projective complex irreducible holomorphic symplectic manifolds of generalized Kummer type and complex dimension \(2n\), where \(n \geq 2\). Then there exists an exact equivalence
  \[D^{b}(X) \cong D^{b}(X').\]
\end{theorem}

\begin{proof}
  We follow the argument of \cite[Section~2.6]{MSYZ2025}. By \Cref{cor:adjacent-birational-models}, it suffices to treat the case in which the K\"ahler cones of \(X\) and \(X'\), identified by the birational map, share a wall supported by \(w^{\perp}\) for a primitive algebraic MBM class \(w\). The isomorphic case is immediate, and the general case follows by composing equivalences along a finite chain of adjacent birational models.

  Choose \(D,g\) and the parallel-transport isometry \(\rho\) as above, and set
  \[\Psi=\rho\circ\psi_E:H^2(K(M),\QQ)\xrightarrow{\sim}H^2(X,\QQ).\]
  Fix the initial markings \(\eta_M,\eta_A\), so that \(\psi=\eta_A\psi_E\eta_M^{-1}\), and let \(\cM_\psi^0\) be the component containing the initial pair. Give \(X\) the marking \(\eta_X=\eta_A\rho^{-1}\), and give \(X'\) the marking \(\eta_{X'}\) induced by the birational map. Since these identifications are parallel transport, both marked targets lie in the target component of \(\cM_\psi^0\). By the surjectivity in \Cref{lem:pair-paths}, we can complete them to two quadruples
  \[
    (Y,\eta_Y,X,\eta_X),\qquad
    (Y',\eta_{Y'},X',\eta_{X'})\in\cM_\psi^0.
  \]
  The marked targets are inseparable by \Cref{thm:torelli}. The marked sources \((Y, \eta_{Y}), (Y', \eta_{Y'})\) also have the same period and lie in the same connected component, so they too are inseparable. Using their markings to identify cohomology spaces with \(H^2(K(M),\QQ)\), their common period is
  \[l_S=\Psi_{\C}^{-1}(H^{2,0}(X)).\]
  With the target cohomology identified birationally, both Hodge isometries are represented by \(\Psi\).

  By \Cref{prop:no-wall}, \(\Psi^{-1}(w)\) is not proportional to an integral MBM class. Hence \(\Psi^{-1}(w^{\perp})\) is not a source wall hyperplane. Local finiteness of the MBM hyperplanes, as in \Cref{lem:mbm-local-finiteness}, allows us to choose a point \(x_0\) in the relative interior of the common target wall, lying on no other target wall, such that \(\Psi^{-1}(x_0)\) lies in the interior of a source chamber \(\mathcal C\). A small neighborhood of \(\Psi^{-1}(x_0)\) stays in \(\mathcal C\). We can therefore choose
  \[
    x_+\in\cK_X,\qquad x_-\in\cK_{X'},\qquad
    \Psi^{-1}(x_+),\Psi^{-1}(x_-)\in\mathcal C.
  \]

  By \Cref{thm:mbm-chambers}, there is a birational model \(Y''\) of \(Y\) with a marking \(\eta_{Y''}\) in the same source component and with the same period, whose K\"ahler cone represents \(\mathcal C\). Writing
  \[\tau=\eta_M^{-1}\eta_{Y''}:H^2(Y'',\Z)\xrightarrow{\sim}H^2(K(M),\Z),\]
  we have \(\tau(\cK_{Y''})=\mathcal C\). Thus \(\tau^{-1}\Psi^{-1}\) sends \(x_+\) and \(x_-\) to K\"ahler classes on the same \(Y''\), and both quadruples with source \((Y'',\eta_{Y''})\) are K\"ahler compatible. The marked point \((Y'',\eta_{Y''})\) is inseparable from both \((Y,\eta_Y)\) and \((Y',\eta_{Y'})\) by \Cref{thm:torelli}. Applying \Cref{lem:inseparable-pairs} to the two quadruples constructed above gives
  \[
    (Y'',\eta_{Y''},X,\eta_X),\qquad
    (Y'',\eta_{Y''},X',\eta_{X'})\in\cM_\psi^0.
  \]

  The rational Hodge isometry \(\tau^{-1}\Psi^{-1}\) pulls an ample class of \(X\) back to a rational \((1,1)\) class of positive square on \(Y''\). Clearing denominators and applying \Cref{thm:projectivity} shows that \(Y''\) is projective. By \Cref{prop:ext-continue}, the two quadruples therefore yield twisted derived equivalences. Since both continuations transport \(B_M\), their source Brauer classes coincide to \(\tau^{-1}(B_M)\). Denote this class by
  \[\beta=[\tau^{-1}(B_M)]_{Y''}.\]
  On the targets, \(\rho(B_A)=D/(cg)\) is algebraic on both \(X\) and \(X'\), so both target Brauer classes vanish. Consequently, the two equivalences compose to give
  \[D^b(X)\cong D^b(Y'',\beta)\cong D^b(X'),\]
  as required.
\end{proof}

\bibliographystyle{alpha}
\bibliography{lib}

\end{document}